\documentclass[11pt]{amsart}

\usepackage{fancybox,amscd,eucal,enumerate}
\usepackage{amsmath,amsthm,amssymb,amsfonts}
\usepackage[all]{xy}
\usepackage{array,tabularx}
\usepackage{hyperref}

\usepackage{paralist}

\newtheorem{theorem}{Theorem}[section]
\newtheorem{lemma}[theorem]{Lemma}
\newtheorem{proposition}[theorem]{Proposition}
\newtheorem{corollary}[theorem]{Corollary}

\theoremstyle{definition}

\newtheorem{definition}[theorem]{Definition}

\newtheorem{properties}[theorem]{Properties}
\newtheorem{facts}[theorem]{Facts}
\newtheorem{claim}[theorem]{Claim}

\newtheorem{remark}[theorem]{Remark}
\newtheorem{notation}[theorem]{Notation}

\theoremstyle{remark}

\newenvironment{penumerate}[1][]{\parindent=0 in \begin{asparaenum}[#1]}{\end{asparaenum}}

\newcommand{\wPsi}{\widehat{\Psi}}
\newcommand{\oPsi}{\overline{\Psi}}

\newcommand{\field}[1]{\ensuremath{\mathbb{#1}}}
\newcommand{\bba}{\field{A}}
\newcommand{\bbc}{\field{C}}
\newcommand{\bbf}{\field{F}}
\newcommand{\bbh}{\field{H}}

\newcommand{\bbp}{\field{P}}
\newcommand{\bbq}{\field{Q}}
\newcommand{\bbr}{\field{R}}

\newcommand{\bbz}{\field{Z}}

\newcommand{\uc}{\underline{\mbc}}

\newcommand{\uz}{\underline{\bbz}}

\newcommand{\uA}{\underline{A}}
\newcommand{\uM}{\underline{M}}

\newcommand{\oJ}{\bar{J}}

\newcommand{\oi}{\bar{\imath}}

\newcommand{\ow}{\bar{w}}

\newcommand{\wJ}{\widehat{J}}

\newcommand{\cala}{\mathcal{A}}
\newcommand{\calb}{\mathcal{B}}
\newcommand{\calc}{\mathcal{C}}

\newcommand{\cale}{\mathcal{E}}

\newcommand{\calh}{\mathcal{H}}
\newcommand{\cali}{\mathcal{I}}

\newcommand{\calm}{\mathcal{M}}
\newcommand{\caln}{\mathcal{N}}
\newcommand{\calo}{\mathcal{O}}

\newcommand{\calq}{\mathcal{Q}}
\newcommand{\calr}{\mathcal{R}}
\newcommand{\cals}{\mathcal{S}}
\newcommand{\calt}{\mathcal{T}}
\newcommand{\calv}{\mathcal{V}}

\newcommand{\bone}{\mathbf{1}}

\newcommand{\mbc}{\mathbf{c}}

\newcommand{\mbf}{\mathbf{f}}

\newcommand{\mbh}{\mathbf{h}}

\newcommand{\mbq}{\mathbf{q}}

\newcommand{\mbu}{\mathbf{u}}

\newcommand{\mbx}{\mathbf{x}}
\newcommand{\mbz}{\mathbf{z}}

\newcommand{\scx}{\text{\sc x}}
\newcommand{\scy}{\text{\sc y}}

\newcommand{\sct}{\text{\sc t}}
\newcommand{\scX}{\text{\sc X}}

\newcommand{\la}{\langle}
\newcommand{\ra}{\rangle}

\newcommand{\dl}{[ \! [}
\newcommand{\dr}{] \! ]}

\newcommand{\cxs}{\Sigma \!\!\!\! /\,}

\newcommand{\equdef}{:=}

\newcommand{\id}{\operatorname{id}}                       % Identity map

\newcommand{\abor}[3]{H_{\text{bor}}^{#1,#2}(#3,\bbz)}  % Assoc. Borel cohomology - Z

\newcommand{\aborM}[3]{H_{\text{Bor}}^{#1,#2}(#3,\uM)}  % Assoc. Borel cohomology - M

\newcommand{\symg}[1]{{{\mathfrak S}_{#1}}}

\newcommand{\dual}[1]{{#1}^\vee}

\newcommand{\ve}{\varepsilon}

\newcommand{\si}{s_{\infty}}
\newcommand{\sih}{s_{\infty *}}
\newcommand{\sic}{s^*_\infty}
\newcommand{\jd}{{\rm j}_{_{\dagger} } }

 \date{April 1st, 2004}
 \title[Equivariant cohomology of real quadrics]{Bigraded Equivariant Cohomology \\ of Real Quadrics}
 \author[dos Santos]{Pedro F. dos Santos}
 \address{Departamento de Matem\'atica, Instituto Superior
 T\'ecnico, Portugal}
 \email{pedfs@math.ist.utl.pt}
 \author[Lima-Filho]{Paulo Lima-Filho}
 \address{Department of
 Mathematics, Texas A{\&}M University, USA}
 \email{plfilho@math.tamu.edu}
 \thanks{The first author was supported in part by  FCT (Portugal) through program POCTI}


\begin{document}

\begin{abstract}
We give a complete description of  the bigraded Bredon cohomology
ring of smooth projective real quadrics, with coefficients in the
constant Mackey functor $ \mathbf{Z} $.  These invariants are
closely related to the integral  motivic cohomology ring, which is
not known for these varieties. Some of the results and techniques
introduced can be applied to other geometrically cellular real
varieties.
\end{abstract}

 \maketitle

 \tableofcontents
%\vfill\eject

%%%%%%%%%%%%%%%%%%%%%%%%%%%%%%%%%%%%%%%%%%%%%%%%%%%%%%%%%%%%%%%%%%%%%%%%%%%%%%%%%%

% \input{intro.tex}
% \vfill\eject

\section{Introduction}
\label{sec:intro}

This paper provides a complete description the bigraded
equivariant Bredon cohomology ring of smooth projective real
quadrics, with coefficients in the constant MacKey functor $\uz$.
From the geometric point of view,  the main motivation for this
work is that these cohomology rings can detect algebraic geometric
invariants that are invisible to ordinary cohomology. From the
topological point of view, the results and techniques presented
here can be extended to a vast class of spaces and constitute a
rather non-trivial family of examples in equivariant cohomology.

Let $\symg{2}$ denote the Galois group $Gal(\bbc/\bbr)$ and let
$RO(\symg{2})=\bbz\cdot \bone \oplus \bbz\cdot \sigma$ be its real
representation ring, generated by the trivial representation
$\bone $ and the \emph{sign} representation $\sigma$. The Bredon
cohomology of a $\symg{2}$-space $Y$ with coefficients in $\uz$ is
an $RO(\symg{2})$-graded ring, written additively as
\begin{equation}
\label{eq:bredonY}
H^*(Y;\uz) = \oplus_{\alpha \in RO(\symg{2})}\ H^\alpha(Y;\uz),
\end{equation}
whose multiplication sends $H^\alpha(Y;\uz)\otimes H^\beta(Y;\uz)$
into $H^{\alpha + \beta}(Y;\uz)$, where $\alpha+\beta$ denotes
addition in the representation ring. We adopt the \emph{motivic
notation} under which we denote $H^{(r-s)\bone + s \sigma}(Y;\uz)$
by $H^{r,s}(Y;\uz)$, for integers $r$ and $s.$

When $X$ is a real algebraic variety, $\symg{2}$ acts via complex
conjugation on the space $X(\bbc)$ of complex points of $X$
endowed with the analytic topology. In this context, an important
property of this bigraded Bredon cohomology is the existence of
natural \emph{cycle maps}
\begin{equation}
\label{eq:cycle}
\gamma \colon H_{\calm}^{r,s}(X;\bbz) \to H^{r,s}(X(\bbc);\uz),
\end{equation}
from the \emph{motivic cohomology} $H_\calm^{r,s}(X;\bbz)$ of $X$
to the Bredon cohomology of $X(\bbc)$, cf. \cite{Dug&Isa-hyper},
which assemble into a bigraded ring homomorphism
\begin{equation}
\label{eq:ring_cycle}
\gamma \colon H_{\calm}^{*,*}(X;\bbz) \to H^{*,*}(X(\bbc);\uz).
\end{equation}
In particular, this map can become a useful tool to detect
non-trivial motivic cohomology classes of real varieties. The
primary example in this paper is the $n$-dimensional real quadric
$\scX_\mbq$ associated to a non-degenerate real quadratic form
$\mbq$ of rank $n+2$ and Witt index $s$. Surprisingly, even the
classical Chow ring $CH^*(\scX_\mbq) \cong \oplus_{p\geq 0}\
H_\calm^{2p,p}(\scX_\mbq;\bbz)$  of $\scX_\mbq$ as a real variety
seems to be unknown.

Let $\calb$ denote the Bredon cohomology ring $\calb^{*,*}\equdef
H^{*,*}(pt;\uz)$ of a point;  see
Section~\ref{sub-sec:cohomology-point}. Its subring $\calb_+ =
\oplus_{r\geq 0, s\geq 0}\ \calb^{r,s}$ generated by the
non-negatively graded elements is isomorphic to a polynomial ring
$\aleph[\tau]$, where $\aleph = \bbz[\ve]$ with $2\ve=0$. Under
this identification, $\ve$ is the generator of $\calb^{1,1}\cong
\bbz/2\bbz\ $ and $\tau$ is the generator of $\calb^{0,2}\cong
\bbz.$

Given indeterminates $h, y$ and $\scx$, let $\la p_1,\ldots,p_k
\ra $ denote the ideal generated by elements $p_1,\ldots,p_k$ in
the polynomial ring $\calb[h,\scx,y].$ Now, given $0\leq s \in
\bbz$, define a subring $\calb_s[h,\scx,y]$ of $\calb[h,\scx,y]$
by
\begin{equation}
\label{eq:def_Bs}
\calb_s[h,\scx,y] \equdef \calb[h]+ \la h^s \ra .
\end{equation}
If $p_1,\ldots,p_k$ are elements in $\calb_s[h,\scx,y]$, define
\begin{equation}
\label{eq:id-notation}
[p_1,\ldots,p_k] \equdef \la p_1,\ldots,p_k \ra \cap \calb_s[h,
\scx,y]
\end{equation}
and denote $\calb_s[h,y] \equdef \calb_s[h,\scx,y] \cap
\calb[h,y].$

Setting $\deg{h} = (2,1)$ define for each nonnegative integer $m$
the following bihomogeneous polynomial of degree $(2m+1,2m+1)$:
\begin{equation}
\label{eq:polys-intro}
\mbf_m \equdef \sum_{a+2b=m}\ \binom{a+b}{b} \ve^{2a+1} \tau^b
h^{2b} \ \in \ \calb_+[h],
\end{equation}
whose properties are studied in Section \ref{sec:alg}. These
polynomials are used in the description of the main result below,
where $\Lambda(\eta)$ denotes the exterior algebra over $\bbz$ on
a single element $\eta$.
\medskip

\noindent{\bf Theorem A.}\ {\it Let $\scX_\mbq$ be the
$n$-dimensional real quadric associated to a non-degenerate
quadratic form $\mbq$ of rank $n+2$ and Witt index $s, $ with
$n\geq 2s$. If $n=2m-\delta,$ with $\delta \in \{ 0, 1 \}$ and
$\deg{\eta} = (2(n-s+1),n-s+1)$, $\deg{h} = (2,1)$,
$\deg{\scx}=(n,-1)$ and $\deg{y}=(0,-2)$, one has a ring
isomorphism
$$ H^{*,*}(\scX_\mbq(\bbc);\uz) \cong
\calb_s[h,\scx,y]\otimes \Lambda(\eta)/ \cali_{n,s},
$$
where ideal $\cali_{n,s}$ can be written as
$$\ \cali_{n,s} = [h^s]\cdot \tilde{J}_{n-2s} \otimes \Lambda(\eta)\ +\ [h^s]\otimes
\la \eta \ra  \ + \ \la h^{n-s+1}\otimes 1 - 2(1\otimes \eta) \ra
,\
$$
where $\tilde{J}_{n-2s} = [g_1,g_2, g_3,g_4,g_5] \subset
\calb[h,\scx,y]$ is the ideal generated by the elements
\begin{align*}
g_1 & = \mbf_{m-s}, \\ g_2 & = \ve^{1-\delta} \tau^{m-s} \scx -
h^{1-\delta} \mbf_{m-s-1}, \\ g_3 & = h \scx, \\
g_4& = h^{2(m-s)} - \delta \{ (-1)^{m-s} \tau^{m-s+1} \scx^2 \} ,
\ \ \text{ and } \\ g_5 & = \tau y - 1 .
\end{align*}  }
\smallskip

\noindent The elements $h, \scx$ and $\eta$ have explicit
geometric origin, described in detail in Sections \ref{sec:aniso}
and \ref{sec:isotropic}.

In Section \ref{sec:isotropic} we discuss this ring structure from
another perspective more suitable for explicit calculations. In
some particular cases, such as \emph{Pfister quadrics} or general
\emph{anisotropic quadrics}, this presentation acquires a simpler
form. Let us first introduce the \emph{associated Borel cohomology
ring} $\cala := H^{*,*}(E\symg{2};\bbz)$, where $E\symg{2}$ is the
classifying space of $\symg{2}.$ This ring is obtained from
$\calb$ by inverting the element $\tau,$ i.e. $\cala\cong
\calb[y]/\la y\tau-1 \ra $, and is isomorphic to
$\aleph[\tau,\tau^{-1}];$ see Remark \ref{rem:fg}(ii).
\medskip

\noindent{\bf Theorem B.\ }{\it Let $\mbq$ be an anisotropic real
quadratic form of rank $n$. Then:
\begin{enumerate}
\item For $n=2m-1$ one has an isomorphism of $\calb$-algebras:
$$
H^{*,*}(\scX_\mbq(\bbc);\uz) \cong \cala[h]/\, I_{2m-1},
$$
where $I_{2m-1}\ = \ \la\, \mbf_m,\ h\mbf_{m-1},\ h^{2m}\, \ra$,
with $\deg{h} = (2,1)$.
\item For $n=2m$   one has a ring isomorphism:
$$
H^{*,*}(\scX_\mbq(\bbc);\uz) \cong \cala[h,\scx]/\, J_{2m},
$$
where $J_{2m}\ = \ \la\, \mbf_m ;\ \ve \tau^m \scx - h
\mbf_{m-1};\ h\scx;\ h^{2m}-(-1)^m\tau^{m+1}\scx^2\, \ra$, with
$\deg{h} = (2,1)$ and $\deg{\scx}=(2m,-1)$.
\end{enumerate}
}
\smallskip

\noindent{\bf Corollary C.\ }{\it If $\scX$ is a Pfister quadric
of dimension $2^{r+1}-2$, then
$$
H^{*,*}(\scX_\mbq(\bbc);\uz) \cong \cala[h,\scx]/\, J_{2^{r+1}-2},
$$
where $J_{2^{r+1}-2}$ is the ideal
$$
\la\, \ve^{2^r-1} ;\ \ve \tau^{2^r-1} \scx - h
\ve\sum_{j=0}^{r-2}\ (\ve^4)^{2^r-2^j} \tau^{2^j-1}(h^2)^{2^j-1};\
h\scx;\ h^{2^{r+1}-2} + \tau^{2^r}\scx^2\, \ra .
$$}
\smallskip

The paper is organized as follows. Section \ref{sec:top-prelim}
provides the necessary topological background, for the reader less
familiar with Bredon cohomology and its associated Borel
cohomology. In Section \ref{sec:geocell} we introduce spectral
sequences
$$ E_*^{i,j}(p)
\Rightarrow \abor{i+j}{p}{Y}
$$
converging to the associated Borel cohomology - see Definition
\ref{def:gcell} - of a $\symg{2}$-space $Y.$ These sequences have
a ``tri-graded'' multiplicative structure and acquire a
particularly interesting form when $Y=X(\bbc)$ where $X$ is a
\emph{geometrically cellular} real variety; cf. Definition
\ref{def:gcell}. Given a real variety $X$, let
$\calc\calh^k(X_\bbc)$ denote the Chow group $CH^k(X_\bbc)$ of
codimension-$k$ algebraic cycles modulo  rational equivalence,
seen as a $\symg{2}$-module under the action of the Galois group
on cycles. Also, for $m \in \bbz$, let $\bbz(m)$ denote the
$\symg{2}$-submodule $(2\pi i)^m\bbz \subset \bbc$, and for a
$\symg{2}$-module $\caln$ let $H^*(\symg{2};\caln)$ denote group
cohomology with coefficients in $\caln$.
\smallskip

\noindent{\bf Proposition D.}\ {\it Let $X$ be a
  geometrically cellular real variety and let $X_\bbc$
  be the complex variety obtained by base-change. Then there is a family
  of spectral sequences $\{ E_r^{*,*}(p) , d_r \}$ converging to
{\rm $\abor{*}{p}{X(\bbc)}$}, with
  $$
  E_r^{i,j}(p) = \begin{cases}
  H^i\left(
  \symg{2};\calc\calh^k(X_\bbc)\otimes \bbz(q-k)
  \right) &, \text{ if } j=2k
  \text{ is even }
  \\ 0 &, \text{ if } j \text{ is odd}.  \end{cases}
  $$
}
\medskip

It must be noted that B. Kahn has constructed in \cite{Kahn-cell}
a spectral sequence converging within a range to Lichtenbaum's
\'etale motivic cohomology of a real variety whose $E_2$-term
coincides with ours. A comparison between these two spectral
sequences is yet to be done.

Our strategy is to start with anisotropic quadrics in Section
\ref{sec:aniso}. In this case Bredon cohomology coincides with its
associated Borel theory and one can use the spectral sequences
above to aid the computation. After the identification of the
appropriate generators $h$ and $\scx$ in  Definition
\ref{def:Psi}, we use the various algebraic results from Section
\ref{sec:alg} to completely determine the relations defining the
ideal $J_{n}$, proving Theorem B above.

In Section \ref{sec:isotropic} we deal with arbitrary isotropic
quadrics. In this case, the additive structure of the cohomology
has a classical decomposition which, in current terminology,
follows directly from the decomposition of the \emph{motive} of
the quadric. Roughly speaking, if $\mbq = \mbq' + \mbh$ is a
quadratic form of rank $n+2$, where $\mbq'$ is anisotropic and
$\mbh$ is hyperbolic of rank $s$, then
$H^{*,*}(\scX_\mbq(\bbc);\uz) \cong \ $ is isomorphic to
$$
\left( \oplus_{j=0}^{s-1}\ \calb\cdot \mbh^j \right) \bigoplus \
H^{*-2s,*-s}(\scX_{\mbq'};\uz)\ \bigoplus \left(
\oplus_{j=0}^{s-1}\ \calb \cdot \eta \mbh^j \right),
$$
where $\mbh$ is the first Chern class in Bredon cohomology of the
hyperplane bundle and $\eta$ is the Poincar\'e dual of a maximal
real linear subspace of dimension $s-1$ contained in $\scX_\mbq.$
One should contrast this decomposition with Theorem $A.$  The
multiplicative structure involves a careful study of the maps
involved in this motivic splitting, along with the relationship
between the coefficient rings $\cala$ and $\calb$.
\smallskip

%%%%%%%%%%%%%%%%%%%%%%%%%%%%%%%%%%%%%%%%%%%%%%%%%%%%%%%%%%%%%%%%%%%%%%%%%%%%%%%%%%

%\input{newpreliminaries.tex}
%\vfill\eject

\section{Background}
\label{sec:top-prelim}

This section contains a brief summary of the main properties of
$RO(\symg{2})$-graded Bredon cohomology.

\begin{definition}
\label{def:modules}
Given $q\in \bbz$, define $\bbz(q):= (2\pi i)^q\bbz \subset \bbc $
with the $\symg{2}$-module structure induced by complex
conjugation. If $M$ is a $\symg{2}$-module, denote $M(q)\equdef
M\otimes_\bbz \bbz(q)$, and let $\underline{M}(q)$ be the
associated Mackey functor. For simplicity, write $\uM \equdef
\uM(0).$ We denote by $\bbz[\xi,\xi^{-1}]$  the
$\bbz[\symg{2}]$-subalgebra of $\bbc$ generated by $\xi:=2\pi i$.
It has a natural graded ring structure defined by setting $\deg\xi
=1$. As a $\bbz[\symg{2}]$-module, we have $ \bbz[\xi,\xi^{-1}]
\cong \oplus_{q\in \bbz} \bbz(q)$. If $M$  is a
$\bbz[\symg{2}]$-module, let $M[\xi,\xi^{-1}]$ denote the
$\bbz[\symg{2}]$-module $M\otimes_\bbz \bbz[\xi,\xi^{-1}]$.
\end{definition}

Given a $\symg{2}$-module $M$ and a $\symg{2}$-space $X$, the
$RO(\symg{2})$-graded Bredon cohomology of $X$ with coefficients
in $\uM$  is an $RO(\symg{2})$-graded abelian group
\[
H^*(X;\uM) = \oplus_{\alpha \in RO(\symg{2})}\ H^\alpha(X;\uM).
\]
We adopt the \emph{motivic notation} under which $H^{(r-s)\bone +
s \boldsymbol{\sigma}}(X;\uz)$ is denoted $H^{r,s}(X;\uz)$, where
$r,s\in\bbz$, $\bone$ denotes the trivial representation of
dimension $1$ and $\boldsymbol{\sigma}$ denotes the  \emph{sign
representation}.

\begin{properties}
\label{prop:bredon}
 Fix a $\symg{2}$-module $M$.
\begin{penumerate}[\bf i.]
\item\label{prop:bredon:i} There is a \emph{forgetful functor}
$\varphi \colon H^{p,q}(X;\underline{M}) \to H^p(X;M(q))$ to
ordinary singular cohomology. This maps factors as
$$H^{p,q}(X;\uM) \to
\calh^p(X;M(q))^{\symg{2}} \hookrightarrow H^p(X;M(q)),$$ where
$\calh^p(X;M(q))^{\symg{2}}$ denotes the invariants of
$H^p(X;M(q))$ considered as a $\symg{2}$-module via the
simultaneous action of $\bbz[\symg{2}]$ on both $X$ and $M.$
\item\label{prop:bredon:ii} There is a {\em transfer functor}
$\tau\colon  H^p(X;M)\to H^{p,q}(X;\underline{M})$ such that the
composite $\tau\circ\varphi\colon H^{p,q}(X;\underline{M})\to
H^{p,q}(X;\underline{M})$ is multiplication by $2$.
\item\label{prop:bredon:iii}  If $A$ is a commutative
$\bbz[\symg{2}]$-algebra, the multiplication
$A\otimes_{\bbz[\symg{2}]} A\to A$ induces a structure of
$(\bbz\times \bbz)$-graded ring on
$$H^{*,*}(X;\underline{A})\equdef \bigoplus_{(p,q)\in \bbz\times
\bbz}\ H^{p,q}(X;\underline{A}).$$ In this case,
$H^*(X;A[\xi,\xi^{-1}])$ has a natural struture of
$\bbz[\symg{2}]$-algebra  and the forgetful functor becomes a map
of (bigraded) $\bbz[\symg{2}]$-algebras
$$
\varphi\colon H^{*,*}(X;\underline{A}) \to H^*(X;A[\xi,\xi^{-1}]),
$$
whose image lies in the invariant subring
$\calh^*(X;A[\xi,\xi^{-1}])^{\symg{2}}.$

\item\label{prop:bredon:iv} There is a natural isomorphism
$H^{p,0}(X;\uM) \cong H^p(X/\symg{2};M).$

\item\label{prop:bredon:v}
If $\symg{2}$ acts freely on $X$ and $2M=0$ there is a natural
isomorphism $H^{p,q}(X;\uM) \cong H^p(X/\symg{2};M(q))$.

\item\label{prop:bredon:vi}
There is a natural isomorphism $H^{p,q}(X\times\symg{2};\uM)\cong
H^p(X;M)$. Under this isomorphism the forgetful functor $\varphi$
is identified with the map $pr_1^*\colon H^{p,q}(X;\uM)\to
H^{p,q}(X\times\symg{2};\uM)$ induced by the projection
$pr_1:X\times\symg{2}\to X$.

\end{penumerate}
\end{properties}

\begin{definition}
The \emph{associated Borel cohomology} to $H^{*,*}(-;\uM)$ is
defined as
$$
\aborM{p}{q}{X} \equdef H^{p,q}(X\times E\symg{2};\uM) ;
$$ cf. \cite[p. 35]{may}.
\end{definition}

\begin{properties}
\label{prop:borel}
Fix an abelian group $M$.

\begin{penumerate}[\bf i.]
\item\label{prop:borel:i}  The projection   $pr_1 \colon X\times
E\symg{2}\to X$ induces a natural transformation $$pr_1^* \colon
H^{p,q}(X;\uM) \to \aborM{p}{q}{X},$$  which is a ring
homomorphism whenever $M=A$ is a $\bbz[\symg{2}]$-algebra.

\item\label{prop:borel:ii}  When $\symg{2}$ acts freely on $X$,  $pr_1^*$ is an isomorphism.
In particular, if $X$ is a finite $CW$-complex of dimension $m$
then $H^{r,s}(X;\uM)=0,$ for all $r>m.$ Furthermore, whenever
$M=A$ is a $\bbz[\symg{2}]$-algebra then
$H^{*,*}(X;\underline{A})$ becomes an algebra over
$H^{*,*}_\text{Bor}(pt;\underline{A}).$

\item\label{prop:borel:iii}  There is a forgetful functor
$\varphi_\text{Bor} \colon H^{p,q}_\text{Bor}(X;\uM) \to
H^p(X;M(q))$ making the following diagram commute
$$
\xymatrix{ H^{p,q}(X;\uM) \ar[rr]^{pr_1^*} \ar[dr]_{\varphi} & &
 H^{p,q}_\text{Bor}(X;\uM) \ar[dl]^{\varphi_\text{Bor}}\\
 & H^p(X;M(q)) & }
$$

\item\label{prop:borel:vi}  If $A$ is a $\bbz/2[\symg{2}]$-algebra
then for every $\symg{2}$-space $X$ one has a natural isomorphism
of bigraded rings
$$
H^{*,*}_\text{Bor}( X ;\uA ) \cong H^*(E\symg{2}\times_\symg{2}
X;A)\otimes \bbz[\xi,\xi^{-1}],
$$
where $\alpha \in H^r(X;A)$ has bidegree $(r,0)$,
$\deg{\xi}=(0,1)$ and $\deg{\xi^{-1}}=(0,-1).$ In particular, if
$X$ is a free $\symg{2}$-space, this induces an isomorphism of
bigraded rings
$$
H^{*,*}( X ;\uA ) \cong H^*(X / \symg{2}; A )\otimes
\bbz[\xi,\xi^{-1}].
$$
\end{penumerate}
\end{properties}

\subsection{Cohomology ring of a point}
\label{sub-sec:cohomology-point}

Let $\calb \equdef H^{*,*}(pt;\uz)$ denote the coefficient ring of
the cohomology theory $H^{*,*}(-;\uz)$, and let $\cala\equdef
H^{*,*}_\text{Bor}(pt;\uz)$ denote the coefficient ring of the
associated Borel theory $H^{*,*}_\text{Bor}(-;\uz)$. Denote by
$\pi \colon \calb \to \cala$ the natural ring homomorphism induced
by the projection $\pi \colon E\symg{2} \to pt$; cf. Property
\ref{prop:borel}[i].

Consider $\aleph = \bbz[\ve]\equdef \bbz[X]/\la 2 X \ra $, where
$\ve$ is an indeterminate of bidegree $(1,1),$ satisfying
$2\ve=0$. It follows from Properties~\ref{prop:borel} and basic
computations in group cohomology that $\cala$ has the following
bigraded ring structure.
\begin{equation}
\label{eq:cala}
\cala \cong \bbz[\ve,\tau,\tau^{-1}] = \aleph[\tau,\tau^{-1}] ,
\end{equation}
where $\tau$ has bidegree $(0,2).$

In order to describe $\calb,$ first consider indeterminates $\ve,
\ve^{-1}, \tau, \tau^{-1}$ satisfying $\deg{\ve}=(1,1),\
\deg{\ve^{-1}}=(-1,-1),\ \deg{\tau}=(0,2)$ and
$\deg{\tau^{-1}}=(0,-2).$ Henceforth, $\ve$ and $\ve^{-1}$ will
always satisfy $2\ve=0=2\ve^{-1}.$

As an abelian group, $\calb$ can be written as a direct sum
\begin{equation}
\label{eq:M}
\calb \ \equdef \ \bbz[\ve,\tau]\cdot 1 \ \oplus \ \bbz[\tau^{-1}]
\cdot \alpha  \oplus \ \bbz[\ve^{-1},\tau^{-1}] \cdot \theta
\end{equation}
where each summand is a free bigraded module over the
corresponding ring. The bidegrees of the generators  $1$, $\alpha$
and $\theta$ are, respectively, $(0,0)$, $(0,-2)$ and $(0,-3).$

The product structure on $\calb$ is completely determined by the
following relations
\begin{equation}
\label{eq:rels}
\alpha\cdot \tau = 2, \quad \quad \alpha\cdot \theta =\alpha\cdot
\ve = \theta \cdot \tau = \theta \cdot \ve \ =\ 0 ,
\end{equation}
and the bigraded ring homomorphism $\pi \colon \calb \to \cala$ is
determined by
$$
\begin{array}{llll}
\ve & \mapsto \ve    &\quad   \tau  & \mapsto \tau  \\
\tau^{-j}\alpha & \mapsto 2 \tau^{-j-1}  & \quad
\ve^{-j-1}\alpha &\mapsto 0 \\
\ve^{-j}\theta & \mapsto 0, &\quad  \tau^{-j}\theta & \mapsto 0,
\end{array}
$$
for $j\geq 0.$

\begin{remark}
\label{rem:fg}
\begin{penumerate}[\bf i.]
\item
Note that $\calb$ is not finitely generated as a ring, and that
$\calb$ has no homogeneous elements in degrees $(p,q)$ when
$p\cdot q < 0$.
\item
The ring $\cala$  is obtained from $\calb$ by inverting the
element $\tau,$ i.e. $\cala\cong \calb[y]/\la y\tau-1 \ra $.
\end{penumerate}
\end{remark}

\subsection{Bredon cohomology of real algebraic varieties}

Given a real algebraic variety $X$, we denote by $X(\bbc)$ its set
of \emph{complex points}  endowed with the analytic topology. It
is a $\symg{2}$-space under the action of complex conjugation and
we can consider the its Bredon cohomology ring
$H^{*,*}(X(\bbc);\uz)$. It is related to the motivic cohomology
ring by a homomorphism called the \emph{cycle map} $\gamma \colon
H_{\calm}^{r,s}(X;\bbz)$  $\to H^{r,s}(X(\bbc);\uz)$ which factors
the classical  map
$$
cl\colon CH^n(X)=H_{\calm}^{2n,n}(X;\bbz)\to
H^{2n}(X(\bbc);\bbz(n)).
$$

\begin{definition}
\label{def:gcell}
An algebraic variety $X$ defined over a field $k$ is
\emph{cellular} if there is a filtration $X=X_n\supset
X_{n-1}\supset\dotsb\supset X_0\supset X_{-1}=\varnothing$ by
closed subvarieties such that $X_i-X_{i-1}$ is isomorphic to a
disjoint union of affine spaces $\bba^{n_{ij}}.$ Whenever
$\bar{k}$ is an algebraic closure for $k$ and $X_{\bar{k}}$ is
cellular we say that $X$ is \emph{geometrically cellular}.
\end{definition}

For real cellular varieties the Bredon cohomology ring has a
simple description relating it to the \emph{Chow ring}.

\begin{proposition}
Let $X$ be a cellular real variety. Then there is a natural ring
isomorphism $H^{*,*}(X(\bbc);\uz)\cong
CH^*(X)~\otimes_\bbz~\calb$, where the elements of $CH^n(X)$ are
given degree $(2n,n)$.
\end{proposition}
\begin{proof}
From Definition~\ref{def:gcell} it is easy to show that
$H^{*,*}(X(\bbc);\uz)$ is a free $\calb$-module generated by
elements with bidegrees of the form $(2r,r)$, with $r\geq 0$. Set
$E:=\oplus_{n\geq 0} H^{2n,n}(X(\bbc);\uz)$. Since $\oplus_{n\geq
0}\calb^{2n,n}=\calb^{0,0}\cong\bbz$ it follows that the inclusion
$E^{*,*}\subset H^{*,*}(X(\bbc);\uz)$ induces a graded ring
isomorphism $E^{*,*}\otimes_\bbz\calb^{*,*}\cong
H^{*,*}(X(\bbc);\uz).$

From Definition~\ref{def:gcell} and basic computations in motivic
and Bredon cohomology, it follows that for a real variety $X$ the
cycle map $ \gamma \colon CH^*(X)= \oplus_{n\geq 0}
H_{\calm}^{2n,n}(X;\bbz) \to \oplus_{n\geq
0}H^{2n,n}(X(\bbc);\uz), $ is a ring isomorphism. Hence   $E\cong
CH^*(X)$ and the result follows.
\end{proof}

\section{Descent spectral sequence}
\label{sec:geocell}

For each $\bbz[\symg{2}]$-module $M$ and each integer $q$ there is
a spectral sequence
$$
E_2^{r,s}(q) \equdef H^r\left(\symg{2}; \calh^s(X; M(q)) \right)
\Rightarrow H^{r+s,q}_\text{Bor}(X;\uM) ,
$$
where $H^r\left(\symg{2}; \calh^s(X; M(q))\right)$ denotes group
cohomology of $\symg{2}$ with coefficients in the
$\bbz[\symg{2}]$-module $ \calh^s(X; M(q))$. This sequence can be
seen as the spectral sequence for the homotopy groups of a
homotopy limit cf. \cite{Bou&Kan-homlim}.

\label{prop:borel:v}  If $M=A$ is a $\bbz[\symg{2}]$-algebra,
its multiplication gives rise to a pairing of spectral sequences
\begin{equation}
\label{eq:pairing}
E_2^{r,s}(q)\otimes E_2^{r',s'}(q') \to E_2^{r+r',s+s'}(q+q'),
\end{equation}
for every $q,q'\in \bbz$. This makes  $\{ E^{r,s}_2 \equdef
\bigoplus_{q\in \bbz}E^{r,s}_2(q)\}$ into a spectral sequence of
algebras converging to $H^{*,*}_\text{Bor}(X;\uA)$.

\subsection{Spectral sequences for geometrically cellular
varieties} \hfill

Given a smooth projective variety $X$ defined over a field
$k\subset \bbr$, the cycle map
\begin{equation}
\label{eq:cycle_map}
cl \ : \ CH^*(X_\bbc)\ \longrightarrow\ \bigoplus_{j\geq 0} \
\calh^{2j}(X(\bbc);\bbz(j))
\end{equation}
is a homomorphism of $\bbz[\symg{2}]$-algebras from the
intersection ring $CH^*(X_\bbc)$ to the cohomology ring
$\calh^{2*}(X(\bbc);\bbz(*))$.

When $X$ is a smooth geometrically cellular variety, this map is
an isomorphism. Hence, the $E_2$-term of the spectral sequence
described in  Properties \ref{prop:borel}(iv) becomes
\begin{equation}
\label{eq:ss-cell}
E^{r,s}_2(q) \cong
\begin{cases}
H^r\left( \symg{2};\ CH^j(X_\bbc)\otimes \bbz(q-j) \right) &,
\text{ if } s= 2j,\ 0\leq j \leq \dim{X} \\ 0 &, \text{ otherwise.
}
\end{cases}
\end{equation}

\begin{lemma}
\label{lemma:zerodifferentials}
Let $X$ be a geometrically cellular real variety. If $j\not\equiv
3 \mod{4}$ then the differential $d_j$ for the spectral sequence
$E_*^{*,*}(q)$ vanishes.
\end{lemma}
\begin{proof}
This follows immediately from    \eqref{eq:ss-cell},   basic
computations in group cohomology and the freeness of
$CH^j(X_\bbc)$ for geometrically cellular varieties.
\end{proof}

In order to better understand the ring structure in
\eqref{eq:pairing} for a geometrically cellular variety $X,$
denote
\begin{equation}
\label{eq:ringE}
\cale\equdef \bigoplus_{k\in \bbz} \cale^k,\quad \text{where}
\quad \cale^k \equdef \bigoplus_{r,j}E^{r,2j}_2(j-k).
\end{equation}
Then $\cale$ is a graded ring under the pairing described in
Properties \ref{prop:borel}(iii) and, under this structure,
$\cale^0$ is a subring which inherits its bigrading comes from the
group cohomology ring $H^*(\symg{2};CH^*(X_\bbc))$ of $\symg{2}$
with coefficients in the graded $\bbz[\symg{2}]$-algebra
$CH^*(X_\bbc).$ The elements $\tau \in E^{0,0}(2) \subset
\cale^{-2}$ and $\tau^{-1} \in E^{0,0}(-2) \subset \cale^{2}$
coming from the cohomology ring $\cala$, cf. \eqref{eq:cala},
induce inverse isomorphisms of bigraded $\cale^0$-modules $ \tau
\colon \cale^k \to \cale^{k-2}, $ for all $k\in \bbz.$

Define $\cale_0 = \bigoplus_{k\in \bbz} \cale_0^k  ,$   where
$\cale_0^k = \bigoplus_{j \geq 0} E^{0,2j}(j-k). $ In other words,
$\cale_0$ is the graded subring of $\cale$ formed by the first
columns of the spectral sequences $E^{*,*}(*)$. It follows that $
\cale_0 \cong \left\{ CH^*(X_\bbc)\otimes
\bbz[\xi,\xi^{-1}]\right\}^{\symg{2}}, $ and hence $\cale_0$ is a
free abelian group. In particular, $ \cale^0_0 \cong \bigoplus_j
CH^i(X_\bbc)^{\symg{2}} $ consists of the \emph{Galois invariants}
of the Chow ring $CH^*(X_\bbc)$, and $ \cale^1_0 = \bigoplus_j
CH^i(X_\bbc)^{-}, $ consists of the \emph{anti-invariants}, i.e.
those classes $\alpha \in CH^*(X_\bbc)$ for which $\sigma^*\alpha
= -\alpha.$ Multiplication by $\tau=\xi^2$ induces isomorphisms
$\tau \colon \cale^k_0 \to \cale^{k-2}_0.$

\begin{proposition}
\label{prop:facts}
Let $X$ be a geometrically cellular real variety.

\begin{penumerate}[\bf i.]
\item\label{prop:facts:i}  $2\cdot
H^{*,*}(X(\bbc);\uz)_{\text{tor}}=0$

\item\label{prop:facts:ii}  One has isomorphisms of abelian groups
$$
H^{p,q}_{\text{Bor}}(X(\bbc);\uz) \cong \bigoplus_{r+s=p}
E^{r,s}_\infty(q).
$$
In particular, one has an isomorphism of abelian groups
$$
H^{*,*}_{\text{Bor}}(X(\bbc);\uz) \cong \bigoplus_{j\geq 0}
Gr^jH^{*,*}_{\text{Bor}}(X(\bbc);\uz) = F^j/F^{j+1},
$$
where $H^{*,*}_\text{Bor}(X(\bbc);\uz) = F^0\supseteq F^1
\supseteq \cdots $ denotes the filtration associated to the
spectral sequence. Note that the resulting ring structure on
$E^{*,*}_\infty$ does not determine the ring structure on
$H^{*,*}_\text{Bor}(X(\bbc);\uz).$

\item\label{prop:facts:iii} In the filtration above, each $F^j$ is
an ideal in $F^0 $ and  $Gr^0 \equdef F^0/F^1$ is a free abelian
group. Therefore, the natural projection $F^0\to F^0/F^1$ factors
through $F^0/F^0_\text{tor},$ giving the following commutative
diagram of ring epimorphisms:
$$
\xymatrix{ & H_\text{Bor}^{*,*}(X(\bbc);\uz) \ar[dl] \ar[dr] \\
H_\text{Bor}^{*,*}(X(\bbc);\uz)/\text{tor}\ar[rr]_{\cong} & & Gr^0
H_\text{Bor}^{*,*}(X(\bbc);\uz) }.
$$

\item\label{prop:facts:iv}The forgetful map $\varphi_{Bor} \colon
F^0=H^{*,*}_\text{Bor}(X(\bbc);\uz) \to H^*(X(\bbc);\bbz)$ factors
as  $F^0 \to F^0/F^1 \subseteq \{ CH^*(X_\bbc)\otimes
\bbz[\xi,\xi^{-1}]\}^{\symg{2}} \subseteq
\calh^{2*}(X(\bbc);\bbz(*)).$
\end{penumerate}
\end{proposition}
\begin{proof}
\begin{penumerate}[\bf i.]
\item The assertion follows immediately from Property~\ref{prop:bredon}{.\ref{prop:bredon:ii}}
and the fact that $H^*(X(\bbc);\bbz)$ is torsion free.
\item
Since $H^*(X(\bbc);\bbz)$ is free, it follows from
equation~\eqref{eq:ss-cell} and basic computations  in group
cohomology that $2\cdot E_2^{r,s}=0$ for all $r>0$. That is, all
columns of the $E_2$-term but the first are $2$-torsion. Hence the
same is true for $E_\infty$ and so  $F^1\cong \bigoplus_{j\geq
1}F^j/F^{j+1}$ because $F^j/F^{j+1}$ is an $\bbf_2$-vector space,
for $j\geq 1$. As explained above, the first  column of
$E_2^{*,*}$ is free. Hence the same is true for  $E_\infty$. This
implies that $F^0\cong F^1\oplus F^0/F^1$.

\item The first statement is just a consequence of basic properties of multiplicative
spectral sequences. The second statement is a consequence of the
remarks  preceding this proposition. Now,  $H^*(X(\bbc);\bbz)\cong
F^0/F^1\oplus F^1$ and  $F^0/F^1$ is free. As mentioned in the
proof of \ref{prop:facts:ii}, $F^1$ is torsion, hence
$F^1=H^*(X(\bbc);\bbz)_{\text{tor}}$. This implies that the
diagram commutes and the bottom arrow is an isomorphism.

\item Let $i\colon \symg{2}\to E\symg{2}$ denote the inclusion as the zero skeleton.
By the definition of the filtration $F^\bullet$, $F^1$ is the
kernel of $(\id_X\times i)^*\colon H^{*,*}(X\times
E\symg{2};\uz)\to H^{*,*}(X\times \symg{2};\uz)$. The forgetful
map $\varphi_{\text Bor}$ (cf.
Property~\ref{prop:borel}{.\ref{prop:borel:iii}}), is given by
$\varphi_{\text{Bor}}=(pr_{13}^*)^{-1}\circ pr_{12}^*$, where
$pr_{ij}$ denotes the projection onto the $ij$ factor of $X\times
E\symg{2}\times\symg{2}$. It is easy to see that
$\left(\id_X\times i\right)\circ pr_{13}\simeq pr_{12}$  hence
$\left(\id_X\times i\right)^*=\varphi_{\text{Bor}}$, and
$\varphi_{\text Bor}$ factors as  $F^0\to F^0/F^1 \subset
\calh^*(X(\bbc);\bbz(*))$. Since $X$ is geometrically cellular
$\calh^{2*}(X(\bbc);\bbz(*))$ is identified with
$CH^{*}(X_\bbc)\otimes\bbz[\xi,\xi^{-1}]$. Finally,  by
Property~\ref{prop:bredon}.\ref{prop:bredon:i} we have
$\operatorname{im} \varphi_{\text{Bor}}\subset
\calh^*(X(\bbc);\bbz(*))^\symg{2}$. Together these facts  give a
factorization of $\varphi_{\text{Bor}}$ as stated in the
Proposition.
\end{penumerate}
\end{proof}

%\bigskip

\begin{corollary}
\label{cor:det}
Given a geometrically cellular real variety $X$ then a cohomology
class $\alpha \in H^{*,*}_\text{Bor}(X(\bbc); \uz)$ is completely
determined by its image under the forgetful map $\varphi \colon
H^{*,*}_\text{Bor}(X(\bbc); \uz) \to H^{*}(X(\bbc); \bbz) $ and
the reduction of coefficients map $\rho \colon
H^{*,*}_\text{Bor}(X(\bbc); \uz) \to H^{*,*}_\text{Bor}(X(\bbc);
\underline{\bbf_2}) $
\end{corollary}

\begin{corollary}
\label{cor:reduction}
The restriction of $\rho\colon H^{*,*}(X(\bbc);\uz) \to
H^{*,*}(X(\bbc);\underline{\bbf_2})$ to the torsion subring
$H^{*,*}(X(\bbc);\uz)_{\text{tor}}$ is injective.
\end{corollary}

\begin{remark}
\label{lem:geom-cellular}
If $X$ is a \textbf{cellular} variety, then the spectral sequence
\eqref{eq:ss-cell} degenerates at the $E_2$-term. Examples include
projective spaces $\bbp^n$, Grassmannians $Gr_{r,n}$ and split
quadrics $\calq_{n,n}$.
\end{remark}

%%%%%%%%%%%%%%%%%%%%%%%%%%%%%%%%%%%%%%%%%%%%%%%%%%%%%%%%%%%%%%%%%%%%%%%%%%%%%%%%%%

% \input{algebraic.tex}
% \vfill\eject

\section{Algebraic preliminaries}
\label{sec:alg}

This section contains basic algebraic results, which stem from the
following presentation of the cohomology ring of the real
Grassmannian $Gr_{2,n+2}(\bbr)$ of $2$-planes in $\bbr^{n+2}.$
Consider the polynomial ring $\bbf_2[w_1,w_2]$ on indeterminates
$w_1,w_2$ (Stiefel-Whitney classes) of degrees $1$ and $2$,
respectively. Then
\begin{equation}
\label{eq:grass0}
H^*(Gr_{2,n+2}(\bbr);\bbf_2) \cong \bbf_2[w_1,w_2]/\oJ_n,
\end{equation}
where $\oJ_n = \la \bar{f}_{n+1}, w_2\bar{f}_n \ra$ is the ideal
generated by polynomials defined recursively as $\bar{f}_0 =
1,\bar{f}_1 =w_1,$\ and $\bar{f}_{n+1} = w_1\bar{f}_{n} +
w_2\bar{f}_{n-1}$. See \cite{Hil-coh}.

Recall from section~\ref{sub-sec:cohomology-point} that $\ve$ is
an element of degree $(1,1)$ in the ring $\cala$ and that
$\aleph=\bbz[\ve]$. Consider indeterminates $\xi$ and $h$ of
degrees $(0,1)$ and $(2,1),$ respectively, and define polynomials
$F_m \in \aleph[\xi,h]$ in a similar fashion:
\begin{equation}
\label{eq:polys}
F_0 = 1, \quad F_1 = \ve, \quad \text{ and } \quad F_{m+1} = \ve
F_m + (\xi h) F_{m-1}.
\end{equation}
Note that  $F_m $ is homogeneous  of degree $(m,m).$ One can
assemble these polynomials into generating functions
\begin{equation}
\label{eq:generating}
G\equdef \sum_{m\geq 0} F_m \scy^m \quad \text{and} \quad H\equdef
\sum_{k\geq 0} F_k^2 \scy^{2k} \ \in\ \aleph[\xi,h]\dl \scy \dr
\end{equation}
and verify that
$$
G=\frac{1}{1-p(\scy)} \quad \text{and} \quad
H=\frac{1}{1-p(\scy)^2},
$$
where $p(\scy) = \scy\{ \ve + (\xi h)\ \scy\} .$ In particular $G
= \{ 1+p(\scy)\} H$, and these identities give the following.
\begin{lemma}
\label{lem:id}
For each non-negative integer $m$ one has:
\begin{enumerate}
\item $F_m = \sum_{a+2b=m}\ \binom{a+b}{b}\ve^a \xi^b h^b$.
\item $F_{2m+1} = \ve F_m^2\quad $ and $\quad F_{2m} = F_m^2 + (\xi h)
F_{m-1}^2$. In particular, one can write
$$
F_{2m+1} = \sum_{a+2b=m}\ \binom{a+b}{b}\ve^{2a+1} \xi^{2b} h^{2b}
$$
\end{enumerate}
\end{lemma}

\begin{definition}
\label{def:rings}
Let $\scx_n$ be a variable of degree $(n,-1)$ and let $\calr_n$
denote the bigraded ring $\calr_n\equdef
\aleph[\xi,\xi^{-1},h,\scx_n].$ Note that sending $\tau$ to
$\xi^{2}$ induces an inclusion of bigraded rings $\cala\equdef
\aleph[\tau,\tau^{-1}] \hookrightarrow \aleph[\xi,\xi^{-1}],$ cf.
Definition \ref{def:modules} and \eqref{eq:cala}, which in turn
induces an inclusion
$$
\iota \colon \cala[h,\scx_n]\hookrightarrow \calr_n.
$$
\end{definition}

We now introduce three homomorphisms of bigraded rings. The first
one is a \emph{reduction of coefficients} epimorphism
\begin{equation}
\label{eq:alg_red}
\pi \colon \calr_n \to \bbf_2[\ve,\xi,\xi^{-1},\scx_n,h].
\end{equation}
induced by the quotient map $\aleph \to \aleph/ 2\aleph \cong
\bbf_2[\ve]$. The second one is an isomorphism of bigraded
$\bbf_2$-algebras
\begin{equation}
\label{eq:W}
W \colon \bbf_2[\ve,\xi,\xi^{-1},h,\scx_n] \xrightarrow{\ \cong\ }
\bbf_2[\xi,\xi^{-1},w_1,w_2,\ow_n]
\end{equation}
defined by $\ve \mapsto \xi w_1$, $h \mapsto \xi w_2$ and $\scx_n
\mapsto \xi^{-1} \ow_n$, where $w_1, w_2$ and $\ow_n$ are given
bidegrees  $(1,0), (2,0)$ and $(n,0),$ respectively. The third one
is an epimorphism of bigraded $\bbf_2$-algebras
\begin{equation}
\label{eq:q}
q \colon \bbf_2[\xi,\xi^{-1},w_1,w_2,\ow_n] \to
\bbf_2[\xi,\xi^{-1},w_1,w_2]
\end{equation}
defined by  $\xi\mapsto \xi$, $w_i\mapsto w_i$, $i=1,2$ and
$\ow_n\mapsto \bar{f}_{n}$; where $\bar{f}_n$ is introduced in
\eqref{eq:grass0}.

\begin{remark} For all $n$ and $m$ one can consider the polynomial $F_m$ as an element  in
$\calr_n.$ It follows directly from Lemma~\ref{lem:id}(2) and
basic definitions that:
\label{rem:polys}
\begin{penumerate}[{\bf i.}]
\item $F_{2m+1}$ is the image under $\cala[h,\scx_n] \hookrightarrow \calr_n$ of the polynomial
$\mbf_m \in \cala[h,\scx_n]$, introduced in
\eqref{eq:polys-intro}.
\item $W\circ\pi\,\, (F_n)=\xi^n \bar{f}_{n}$; \ see~\eqref{eq:grass0}.
\end{penumerate}
\end{remark}

\begin{definition}We introduce the  following three ideals:
\label{def:fundamental}

\begin{penumerate}[{\bf i.}]
\item Let $\bar{J}_n$ denote the ideal $\la
\bar{f}_{n+1},w_2 \bar{f}_{n}
\ra\subset\bbf_2[\xi,\xi^{-1},w_1,w_2]$; see~\eqref{eq:grass0}.

\item Write $n=2m-\delta,$ $\delta \in \{ 0, 1\}$, and define
$\wJ_n=~\la \mathbf{g_1}, \mathbf{g_2}, \mathbf{g_3}, \mathbf{g_4}
\ra \subset \calr_n$ as the homogeneous ideal generated by
\begin{enumerate}[]
\item $\mathbf{g_1}:=\ve^{1-\delta}\xi^{2m}\scx_n - h^{1-\delta}F_{2m-1}$
\item $\mathbf{g_2}:=F_{2m+1}$
\item $\mathbf{g_3}:=h \scx_n$
\item $\mathbf{g_4}:=(\xi h)^{2m} -(-1)^m (\xi h)^\delta (\xi^{2m+1-\delta}\scx_n)^2 $
\end{enumerate}

\item Using the inclusion $\iota
\colon\cala[\scx_n,h] \hookrightarrow \calr_n $ define the ideal $
J_n \equdef \cala[\scx_n,h]\cap \wJ_n \ \subset\ \cala[\scx_n,h].
$
\end{penumerate}
\end{definition}

\begin{remark}
\label{rmk:generators-J_n}
Abbreviate $\scx_n$ to $\scx$ for simplicity.
\begin{penumerate}[{\bf (1)}]
\item It is easy to see, from Remark \ref{rem:polys}[i] and simple
degree considerations, that one can present the ideal $J_n \subset
\cala[\scx_n,h]$ as:
$$
J_n =
\begin{cases}
\la \tau^m \scx - \mbf_{m-1},\ \mbf_{m},\ h \scx,\ h^{2m} \ra &,
\text{ if } n =2m-1 \\
\la \ve \tau^m \scx - h \mbf_{m-1},\ \mbf_{m},\ h \scx,\ h^{2m} -
(-1)^m \tau^{m+1} \scx^2  \ra &, \text{ if } n =2m
\end{cases}
$$
\item Note that $\bbf_2[\xi,\xi^{-1},w_1,w_2]/\oJ_n$ is a
presentation of the ring
$$
H^{*}(Gr_{2,n+2}(\bbr);\bbf_2(*))\cong
H^*(Gr_{2,n+2}(\bbr);\bbf_2)\otimes \bbz[\xi,\xi^{-1}] .
$$
\end{penumerate}
\end{remark}

The  commutative diagram below summarizes the relations between
the ideals and maps defined above. The lower vertical maps are the
natural projections to the corresponding quotients, and the
existence of the maps represented by dotted arrows is established
in Proposition~\ref{prop:algebraic}. The composite $q\circ W\circ
\pi$ is denoted by $\Psi.$

{\small
\begin{equation}
\label{eq:alg_diagram}
\xymatrix{
\ve\cala[h,\scx_n] \ar@{^{(}-->}[r]^-{} \ar@{^{(}->}[d]_{} &
\bbf_2[\xi,\xi^{-1},\ve,h,\scx_n] \ar[r]^-{W}_-{\cong}
\ar@{-->}[dr]_-{\overline{\Psi}} &
\bbf_2[\xi,\xi^{-1},w_1,w_2,\bar{w}_n] \ar[d]^{q}      \\
% \ar@/_{4.5pc}/[dd]_{\oq}
\cala[h,\scx_n] \ar@{^{(}->}[r]^-{\iota} \ar@{->>}[d]^{\gamma} &
\calr_n \ar@{->>}[u]_-{\pi} \ar[r]_-{\Psi}  \ar@{->>}[d]_{p}
& \bbf_2[\xi,\xi^{-1},w_1,w_2] \ar@{->>}[d]  \\
\cala[h,\scx_n]/J_n \ar@{^{(}->}[r]^-{\bar\iota}
\ar@/_{1.5pc}/@{-->}[rr]_{\overline\rho} & \calr_n /\widehat{J}_n
\ar@{-->>}[r]^-{r} & \bbf_2[\xi,\xi^{-1},w_1,w_2]/\bar{J}_n
}
\end{equation}
}

In the next section we show that $\cala[h,\scx_n]/J_n$ is a
presentation for the equivariant cohomology ring of
$\scX_\mbq(\bbc)$, where $\mbq$ is anisotropic of rank $n+2$. The
main properties of $\cala[h,\scx_n]/J_n$ and $J_n$ needed for the
proof of this result are stated  in the next proposition.

\begin{proposition}
\label{prop:algebraic}
Using the notation above:

\noindent{\textbf{i.}}\ The composition $\Psi \equdef q\circ
W\circ \pi$ sends $\wJ_n$ into $\oJ_n$. In particular, it descends
to a map $r \colon \calr_n/\wJ_n \to
\bbf_2[\xi,\xi^{-1},w_1,w_2]/\oJ_n$ making the diagram above
commute.

\noindent{\textbf{ii.}}\ Writing $n=2m-\delta$, with $\delta \in
\{ 0, 1 \} $  one has
$$ J_{2m-\delta} + \la \ve \ra = \la\, \ve,\
h^{1-\delta} \scx,\ h^{2m}-(-1)^m\tau^{m+1}\scx^2\, \ra .$$

\noindent{\textbf{iii.}}\ The torsion subgroup $\left(
\cala[\scx_n,h]/J_n\right)_{tor}$ is precisely the ideal $\la
 {\ve} \ra$ generated by the class of $\ve$.

\noindent{\textbf{iv.}}\ The restriction of $\bar\rho \equdef
r\circ \bar\iota \colon \cala[\scx_n,h]/J_n \to
\bbf_2[\xi,\xi^{-1},w_1,w_2]/\oJ_n$ to the torsion subgroup
$\left( \cala[\scx_n,h]/J_n\right)_{tor}$ is injective.

\end{proposition}
\begin{proof}
Using same arguments as in Lemma \ref{lem:id}, one can prove the
following facts concerning the ideal $\oJ_n$:
\begin{facts}
\label{facts:appendix}
\begin{enumerate}[(a)]
\item $w_2^k\oJ_n \subset \oJ_{n+k} \subset \oJ_n$;
\item $\bar{f}_{2n+1}=w_1\bar{f}_{n}^2$;
\item $\bar{f}_{2n}=\bar{f}_{n}^2+w_2\bar{f}_{n-1}^2$.
\end{enumerate}
\end{facts}

Write $n=2m-\delta,$ where $\delta = s(n) \in \{ 0, 1\}$ is the
sign of $n, $ and denote $\scx_n$ by $\scx$.

\noindent{\it \textbf{i.}}\ One just needs to verify the assertion
on the generators of $\widehat{J}_n$.

The generator $\mathbf{g_1}$ is sent to
$\xi^{2m-\delta}\{w_1^{1-\delta} \bar{f}_{2m-\delta}
+w_2^{1-\delta}\bar{f}_{2m-1}\}.$ This expression is $0$ if
$\delta = 1$ and equal to $\xi^{2m} \bar{f}_{2m+1} \in \oJ_n$ when
$\delta = 0$. The generator $\mathbf{g_2}$ is sent to
$\xi^{2m+1}\bar{f}_{2m+1} \in\oJ_{2m}\subset\oJ_n$; cf. Facts
\ref{facts:appendix}(a). The generator $\mathbf{g_3}$ is sent to
$w_2\bar{f}_{n} \in \oJ_n$. Finally, the generator $\mathbf{g_4}$
is sent to $\xi^{4m}\beta_n$ where $\beta_n:=
w_2^{s(n)}(\bar{f}_{n}^2 + w_2^n)$. To show that $\beta_n\in\oJ_n$
we use induction. First, observe that $\beta_1=w_2
\bar{f}_{2}\in\oJ_1$. Also,
$\beta_2=w_1^4=w_1\bar{f}_{2}\in\oJ_2$. Assuming $\beta_k\in\oJ_k$
for $k\leq n$, we can use Facts~\ref{facts:appendix} above to
write
\begin{align*}
\beta_{n+1} & = w_2^{s(n+1)} \left\{
(w_1\bar{f}_{n}+w_2\bar{f}_{n-1})^2
+ w_2^{n+1} \right\} \\
            & = w_2^{s(n+1)} \left( w_1^2\bar{f}_{n}^2 + w_2^2\bar{f}_{n-1}^2 + w_2^{n+1} \right)\\
           & = w_2^{s(n+1)} \left( w_1\bar{f}_{2n+1} + w_2^2\bar{f}_{n-1}^2 + w_2^{n+1}
           \right) \\
            & = w_2^{s(n+1)}w_1\bar{f}_{2n+1}+w_2^2\beta_{n-1}
            \quad \in\quad  \oJ_{2n} + w_2^2\oJ_{n-1}\subset \oJ_{n+1}.
\end{align*}
This completes the proof of {\it \textbf{i.}}.
\smallskip

\noindent{\it \textbf{ii.}}\ This follows directly from inspection
of the generators of $J_n$ noting that $F_{k}\in\la \ve \ra$ for
$k$ odd; cf. Lemma \ref{lem:id}(b).

\noindent{\it \textbf{iii.}}\ It follows directly from (\it
\textbf{ii}) that $\cala[\scx,h]/\la \ve \ra + J_n$ is torsion
free, hence the torsion ideal of $\cala[\scx,h]/J_n$ is $\la \ve
\ra.$

\noindent{\it \textbf{iv.}}\ Let $\psi$, $\bar\psi$ denote the
composites $q\circ W \circ \pi$ and $q\circ{W}$, respectively. By
{\it \textbf{ii.}} it suffices to show that
$$
\psi^{-1}\left( \oJ_n\cap\la w_1\ra\right)\cap\la\ve\ra\subset
\widehat{J}_n ,
$$
because $\psi(\ve)=\xi w_1$. Using Facts \ref{facts:appendix}(b)
and the inductive definition of $ \bar{f}_{n} $ one obtains
$$
\bar{f}_{2m+1}\in\la w_1 \ra \quad \text{and} \quad
\oJ_{2m-\delta} \cap\la w_1\ra=\la \bar{f}_{ 2m+1}, w_2^{2-\delta
}\bar{f}_{2m-1}\ra.
$$
A direct computation gives
\begin{align*}
\psi(F_{2m+1})&        =\xi^{2m+1}\bar{f}_{2m+1},\\
\psi(h^{2-\delta}F_{2m-1})&=\xi^{2m-\delta+1}w_2^{2-\delta}\bar{f}_{2m-1},\\
W(\xi^{2m-\delta+1}\scx-F_{2m-\delta})& =\xi^{2m-\delta}(\bar w_n
-\bar{f}_{n}),
\end{align*}
where  we  denote the image
$\pi(F_n)\in\bbf_2[\ve,\xi,\xi^{-1},h,\scx_n]$ by $F_n$ as well.
Since $q$ is onto and  $\ker q = \la \bar w_n -\bar{f}_{n} \ra$ it
follows that
$$
{\bar\psi\, }^{-1}(\oJ_n\cap\la w_1\ra) = \la F_{2m+1},
h^{2-\delta}F_{2m - 1}, \xi^{2m-\delta + 1}\scx - F_n\ra.
$$
Now, since $F_{2m \pm 1}\in\la \ve \ra$ we conclude that
$$
{\bar\psi\, }^{-1}(\oJ_n\cap\la w_1\ra)\cap\la\ve\ra = \la
F_{2m+1}, h^{2-\delta}F_{2m - 1}, \ve(\xi^{2m-\delta+1}\scx
-F_n)\ra.
$$
Finally, noting that $\ker\pi=2\calr_n$ and $2\ve=0$ we obtain
$$
\psi^{-1}(\oJ_n\cap\la w_1\ra)\cap\la\ve\ra = \la F_{2m+1},
h^{2-\delta}F_{2m - 1}, \ve(\xi^{n+1}\scx
-F_n)\ra\subset\widehat{J}_n.
$$

\end{proof}

It is important to note that in the case of $n$ odd the ring
$\cala[t,\scx_n]/J_n$ simplifies considerably.

\begin{proposition}
\label{prop:odd}
If $n=2m-1$ is odd, define $I_{2m-1}= \cala[h]\cap \wJ_{2m-1}$.
Then  $I_{2m-1}= \la \mbf_{m}, h\mbf_{m-1}, h^{2m} \ra $ and the
inclusion $\cala[h] \underset{i}{\hookrightarrow} \cala[\scx_n,h]$
induces an isomorphism
$$
\oi\colon\cala[h]/I_{2m-1} \to \cala[h,\scx_{2m-1}]/J_{2m-1}.
$$
\end{proposition}
\begin{proof}
This follows directly from Remark \ref{rmk:generators-J_n} and the
fact that $\tau$ is invertible in $\cala$.
\end{proof}

%%%%%%%%%%%%%%%%%%%%%%%%%%%%%%%%%%%%%%%%%%%%%%%%%%%%%%%%%%%%%%%%%%%%%%%%%%%%%%%%%%

%\input{anisotropic.tex}

\section{Anisotropic quadrics}
\label{sec:aniso}

A non-degenerate quadratic form $\mbq$ on a real vector space $V$
of dimension $n+2$ determines a \emph{real smooth $n$-dimensional
quadric hypersurface} $X_{\mbq} \subset \bbp(V)$. Any such
hypersurface is isomorphic to the quadric $\calq_{n,s}\subset
\bbp^{n+1}_\bbr$ given by the equation
\begin{equation*}
\label{eq:quadric}
\mbq_{n,s}(z_0,\ldots,z_{n+1}) = z_0^2+\cdots + z_{n-s+1}^2-
z_{n-s+2}^2 - \cdots - z_{n+1}^2=0,
\end{equation*}
for some $s$ satisfying $0\leq 2s \leq n+2.$ We denote
$\calq_n\equdef \calq_{n,0}   .$

The main theorem in this section computes  the  ring
$H^{*,*}(\calq_n(\bbc);\uz)$, showing that is generated by two
elements over the Borel cohomology of a point. These elements come
from the following geometric constructions.

Let $S_a^{n}$ be the unit sphere in $\bbr^{n+1}$ endowed with the
antipodal $\symg{2}$-action $a : S^{n} \to S^{n}.$ One has a
classical equivariant embedding  $ j : (S^{n}_a ,a) \to
(\calq_{n},\sigma)$ defined by $j(y_0,\ldots,y_{n+1}) = [iy_0:
\cdots : iy_{n+1}:1]$. It is easy to see that $S^{n}_a$ is
$H\uz$-orientable, with fundamental class $[S_a^n]\in
H_{n,n+1}(S_a^n;\uz)$, and hence we obtain Gysin maps
\begin{equation*}
\label{eq:gysin}
j_! \colon H^{r,s}(S^{n}_a;\uz) \to
H^{r+n,s-1}(\calq_{n}(\bbc);\uz) ,
\end{equation*}
since the fundamental class $[\calq_n(\bbc)]$ lives in
$H_{2n,n}(\calq_n(\bbc),\uz)$. Let $\calo(1)$ denote the
hyperplane bundle $\calo_{\calq_n}(1)$, under the embedding
$\calq_n \hookrightarrow \bbp^{n+1}$. The space of complex points
$\calo(1)(\bbc)$ becomes a Real bundle over $\calq_n(\bbc)$, in
the sense of \cite{Ati66}. It follows that $\calo(1)(\bbc)$ has
Chern classes
\begin{equation}
\label{eq:chern}
\uc_i(\calo(1)(\bbc)) \in H^{2i,i}(\calq_n(\bbc);\uz)
\end{equation}
 with values in Bredon cohomology with coefficients in $\uz;$ cf. \cite{dSan-alg}.

Since the action of $\symg{2}$  on $\calq_n(\bbc)$ is free one can
identify $H^{*,*}(\calq_n(\bbc);\uz)$ with
$H_{\text{Bor}}^{*,*}(\calq_n(\bbc);\uz)$, cf.
Properties~\ref{prop:borel}[i], and use its structure of algebra
over $\cala\equdef H^{*,*}_\text{Bor}(pt;\uz)$ in the following
construction.

\begin{definition}
\label{def:Psi}
Let $\Psi_n \colon \cala[h,\scx_n] \to H^{*,*}(\calq_n(\bbc);\uz)$
denote the homomorphism of  bigraded $\cala$-algebras that sends
\begin{align*}
\scx_n\ & \mapsto \ j_! \bone \in H^{n,-1}(\calq_n(\bbc);\uz)
\quad \text{ and } \\
h \ & \mapsto \ \uc_1(\calo(1)(\bbc)) \in
H^{2,1}(\calq_n(\bbc);\uz) .
\end{align*}
\end{definition}

The classes $j_!\bone$ and $\uc_1(\calo(1)(\bbc))$ are completely
determined by their images under the \emph{forgetful map}
$\varphi$ into singular cohomology,  and the \emph{reduction of
coefficients} map $\rho$ into Bredon cohomology with coefficients
in $\underline{\bbf_2};$ cf. Corollary \ref{cor:det}. Hence, it is
useful to identify the compositions
\begin{equation}
\label{eq:forgetful}
\varphi\circ \Psi_n  \colon \cala[h,\scx_n] \xrightarrow{\Psi_n}
H^{*,*}(\calq_n(\bbc);\uz) \xrightarrow{\varphi}
H^*(\calq_n(\bbc);\bbz(*))
\end{equation}
and
\begin{equation}
\label{eq:reduction}
\rho\circ \Psi_n  \colon \cala[h,\scx_n] \xrightarrow{\Psi_n}
H^{*,*}(\calq_n(\bbc);\uz) \xrightarrow{\rho}
H^{*,*}(\calq_n(\bbc);\underline{\bbf_2}).
\end{equation}

Since quadrics are geometrically cellular, the cycle map
\eqref{eq:cycle_map} induces a natural isomorphism
\begin{equation}
\label{eq:iso1}
H^*(\calq_n(\bbc);\bbz(*)) \equdef \oplus_{p,q}
H^p(\calq_n(\bbc);\bbz(q)) \cong CH^*(\calq_{n,\bbc})\otimes
\bbz[\xi,\xi^{-1}],
\end{equation}
where the elements of $CH^j(\calq_{n,\bbc})$ are given bidegree
$(2j,j).$ Furthermore, the Chow groups of complex quadrics are
well-known classical objects, and have the following presentation.

\begin{facts}
\label{prop:chow_ring}
Write $n=2m-\delta,$ with $\delta \in \{ 0, 1 \}$.  Then, the Chow
ring $\ CH^*(\calq_{2m-\delta,\bbc})\ $ is isomorphic to the
graded ring $\ \bbz[h,\phi]/\calc_{2m-\delta}$,\ where $\deg{h}=1,
\ \deg{\phi} = m $, and $\calc_{2m-\delta}$ is the ideal
$$
\calc_{2m-\delta} := \left< h^{1-\delta}(h^m - 2\phi),  \phi^2 -
\frac{1+(-1)^{m } }{2} \  h^{m} \ \phi \right>.
$$
\begin{enumerate}
\item The action
of $\ Gal(\bbc/\bbr)\ $ on $CH^*(\calq_{n,\bbc})$ is determined by
$$
\sigma_*(h)\ =\ h \quad \text{ and } \quad \sigma_*( \phi )\ =\
\frac{1+(-1)^{m}}{2} \ h^{m}\ -\ (-1)^{m} \phi\ .
$$
\item If $i : \calq_{n } \hookrightarrow \calq_{n+1}$ is the
canonical inclusion, then
$$i^* : CH^*(\calq_{n+1,\bbc}) \to
CH^*(\calq_{n,\bbc})$$ is given by
$$
i^*(\phi ) = h^{\frac{1+(-1)^n}{2}}\phi \quad \text{ and } \quad
i^*(h) = h.
$$
\end{enumerate}
\end{facts}

\begin{remark}
\label{rem:ss}
\begin{penumerate}[i.]
\item   Denote the homogeneous coordinates in
$\calq_{2m-\delta}(\bbc)$ by $[\mbz] = [\mbx : x_{2m-\delta} :
x_{2m-\delta+1}],$ where $\mbx = (x_0,\ldots,x_{2m-\delta-1})\in
\bbc^{2m-\delta}.$ Let $L^+, L^- \subset \calq_{2m-\delta,\bbc}$
be the (maximal) linear subspaces given by
{\small
\begin{align*}
\label{eq:planes}
L^+ & = \{ [\mbz] \mid x_0 = (-1)^\delta i x_1,\ x_{2j} =
ix_{2j+1},\ \ \delta \cdot x_{2m-\delta +1} = 0,\ 1\leq j\leq
m-\delta \} \\
L^- & = \{ [\mbz] \mid x_0 = - i x_1,\ \ \ \ x_{2j} = -
ix_{2j+1},\ \ \delta \cdot x_{2m-\delta +1} = 0,\  1\leq j\leq
m-\delta \}.
\end{align*} }
It follows from \cite{Gri&Har-prin} that when $\delta = 0,$ the
cohomology duals of $[L^+]$ and $[L^-]$, respectively, are
distinct classes generating the cohomology group in middle
dimension.

\item It follows from the above description that when $n\not\equiv 0
\mod{4}$, one has $CH^*(\calq_{n,\bbc})^{\symg{2}} =
CH^*(\calq_{n,\bbc})$, and when $n=4m$ one has
\begin{equation}
\label{eq:invs}
CH^*(\calq_{4m,\bbc})^\symg{2} = \bbz[h,h\phi]/\la h^{4m+1} -
2h\phi, (h\phi)^2  \ra
\end{equation}
and
\begin{equation}
\label{eq:anti-invs}
CH^*(\calq_{4m,\bbc})^{-} = \bbz \cdot \{ h^{2m} - 2 \phi \}.
\end{equation}
Note that $h\cdot (h^{2m} - 2\phi) = 0,$ and this determines the
structure of $CH^*(\calq_{4m,\bbc})^{-} $ as a
$CH^*(\calq_{4m,\bbc})^\symg{2}$-module. Therefore, the $E_2$-term
of the spectral sequence \eqref{eq:ss-cell} is given by
{\small $$
\cale =
\begin{cases}
\cala\otimes CH^*(\calq_n) & \text{ if }  n\not\equiv 0 \mod{4} \\
\cala[h,h\phi]/\la h^{4m+1} - 2h\phi, (h\phi)^2  \ra \bigoplus
\cala\cdot \{ h^{2m} - 2\phi \} & \text{ if } n=4m.
\end{cases}
$$ }
\end{penumerate}
\end{remark}

Equation \eqref{eq:iso1} and Facts \ref{prop:chow_ring} provide an
identification
\begin{equation}
\label{eq:chow}
 \calh^*(\calq_n(\bbc);\bbz(*)) \cong
\bbz[\xi,\xi^{-1}][h,\phi]/\calc_n.
\end{equation}
On the other hand, it is well-known that one has a homeomorphism
$$
i_n \colon \calq_n(\bbc)/\symg{2} \to Gr_{2,n+2}(\bbr),
$$
where the latter denotes the Grassmanian of two-planes in
$\bbr^{n+2}$ with its classical topology. In particular, it
follows from \eqref{eq:grass0} and Properties \ref{prop:borel}(ii)
and (vi) that one has  isomorphisms
\begin{align}
\label{eq:quot}
H^{*,*}(\bbq_n(\bbc);\underline{\bbf_2}) & \cong
H^*(Gr_{2,n+2}(\bbr);\bbf_2)\otimes \bbz[\xi,\xi^{-1}] \\ & \cong
\bbf_2[\xi,\xi^{-1},w_1,w_2]/ \oJ_n \notag
\end{align}
where $\oJ_n = \la  \bar{f}_{n+1} ,  w_2 \bar{f}_n \ra$ is
introduced in Definition~\ref{def:fundamental}(i).

\begin{lemma}
\label{lem:classes}
Under the identifications \eqref{eq:chow} and \eqref{eq:quot}, one
has
\begin{enumerate}[i.]
\item \begin{align*}
\rho(\uc_1(\calo(1)(\bbc)) & = \xi w_2 ,  \quad   \quad  \rho(\ve) = \xi w_1  \quad \text{ and } \\
\varphi(\uc_1(\calo(1)(\bbc)) & = h
\end{align*}
\item \begin{align*}
\rho(j_! \bone) &= \ \xi^{-1} \bar{f}_{n} \quad \text{ and } \\
\varphi(j_! \bone) & = \begin{cases}
0 &,\ n=2m-1 \\
\xi^{-m-1}(h^m -2 \phi )&,\ n=2m .
\end{cases}
\end{align*}
\end{enumerate}
\end{lemma}
\begin{proof}
Left to the reader.
\end{proof}

\begin{proposition}
\label{prop:step1}
The map $\Psi_n \colon \cala[h,\scx_n] \to
H^{*,*}(\calq_n(\bbc);\uz)$ factors through  a map $ \oPsi_n
\colon \cala[h,\scx_n]/J_n \to H^{*,*}(\calq_n(\bbc);\uz),$ where
$J_n$ is introduced in Definition \ref{def:fundamental}(iii).
\end{proposition}
\begin{proof}
It follows from Corollary \ref{cor:det} that it suffices to prove
that $J_n \subseteq \ker{\rho\circ \Psi_n}$ and $J_n \subseteq
\ker{\varphi \circ \Psi_n}$. Lemma \ref{lem:classes} shows that
$\rho\circ \Psi_n$ can be factored as the upper triangle in the
following  diagram:

\begin{equation}
\label{eq:comp1}
\xymatrix{ \cala[h,\scx_n] \ar@{^{(}->}[r]^-{\iota} \ar@{->>}[d]
\ar[rrd]_{\rho\circ\Psi_n} & \calr_n =
\bbz[\xi,\xi^{-1},\ve,h,\scx_n] \ar[r]^-{\ q\circ W \circ \pi\ }
& \bbf_2[\xi,\xi^{-1},w_1,w_2] \ar@{->>}[d] \\
\cala[h,\scx_n]/J_n \ar[rr]_{\overline\rho} & &
\bbf_2[\xi,\xi^{-1},w_1,w_2]/\oJ_n, }
\end{equation}
see diagram \eqref{eq:alg_diagram}. It follows from
Proposition~\ref{prop:algebraic}.{\bf i}  that $\rho\circ\Psi_n $
descends to $\bar\rho,$ making the diagram above commute and
proving the first inclusion.

Since $CH^*(\calq_{n,\bbc})$ is torsion-free, Lemma
\ref{lem:classes} shows that the composition $\varphi\circ \Psi_n$
factors as
\begin{equation}
\label{eq:Psi-hat}
\cala[h,\scx_n] \xrightarrow{\wPsi_n} \bbz[\xi,\xi^{-1}][h,\phi]
\to \bbz[\xi,\xi^{-1}][h,\phi]/\calc_n,
\end{equation}
where the last arrow is the quotient map, and $\wPsi_n$ is the map
induced by the composition $\cala \hookrightarrow
\aleph[\xi,\xi^{-1}]\to \bbz[\xi,\xi^{-1}]$ and by sending
$$
h\mapsto h\quad \text{ and } \quad \scx_n \mapsto
\begin{cases} 0
&, n=2m-1
\\ \xi^{-m-1}(h^m - 2\phi) &, n=2m.
\end{cases}
$$

In particular $\la \ve \ra \subseteq \ker (\varphi\circ \Psi_n)$.
As before, write $n=2m-\delta$ with $\delta\in\{0,1\}$ and
abbreviate $\scx_n$ to $\scx$. From the presentation of $J_n$ in
Remark~\ref{rmk:generators-J_n}, it now follows that it suffices
to show that $h^{1-\delta}\scx$ and $h^{2m} -
(-1)^m\tau^{m+1}\scx^2$ are sent to $\calc_{2m-\delta}$ by the
correspondence above.

Computing we get $(\varphi\circ\Psi)(h^{1-\delta}\scx) =
\xi^{-m-1}\left\{ h^{1-\delta}(h^m - 2 \phi) \right\}$ and, using
$\tau = \xi^2$, one can easily verify that
\begin{multline*}
(\varphi\circ\Psi)( h^{2m} - (-1)^m\tau^{m+1}\scx^2) = \\
[1-(-1)^m]h^m(h^m-2\phi) \ - \ 4(-1)^m\left\{ \phi^2 -
\frac{1+(-1)^m}{2} h^m\phi \right\}.
\end{multline*}
Since $m\geq 1$ one concludes that this elements lies in
$\calc_{2m-\delta}.$

 This completes the proof.
\end{proof}

\begin{corollary}
\label{cor:inj1}
The restriction $$\oPsi_{n|_{tor}} \colon
\left(\cala[h,\scx_n]/J_n \right)_{tor} \to
H^{*,*}(\calq_n(\bbc);\uz)$$ of $\oPsi_n$ to the torsion subgroup
is injective.
\end{corollary}
\begin{proof}
It follows from  Corollary \ref{cor:reduction} that $$\rho \colon
H^{*,*}(\calq_n(\bbc);\uz) \to
H^{*,*}(\calq_n(\bbc);\underline{\bbf_2}) $$ is injective on
torsion. Hence, $\oPsi_{n|_{tor}}$ is injective if and only if so
is $\rho \circ \oPsi_{n|_{tor}}$, but the latter statement follows
immediately from diagram \eqref{eq:comp1} and Proposition
\ref{prop:algebraic}.{\bf iv}.
\end{proof}

The following result along with with Proposition \ref{prop:odd},
yield the proof of Theorem B in the introduction.

\begin{theorem}
\label{thm:MAIN}
The map $\oPsi_n \colon \cala[h,\scx_n]/J_n \to
H^{*,*}(\calq_n(\bbc);\uz)$ is an isomorphism of bigraded algebras
over $\cala.$
\end{theorem}

\begin{proof}
The injectivity of $\oPsi_n$ follows from Corollary
\ref{cor:inj1}, Corollary~\ref{cor:det} and the following claim.
\begin{claim}
\label{claim:eq}
$\ker (\varphi\circ \Psi_n) = \la \ve \ra + J_n$
\end{claim}
\begin{proof}
It suffices to prove that $\ker (\varphi\circ \Psi_n) \subseteq
\la \ve \ra + J_n$.

Let us first consider the case $n=2m-1.$ It is clear that any
$p(h,\scx_{2m-1}) \in \cala[h,\scx_{2m-1}]$ is congruent $\mod \la
\ve \ra + J_{2m-1}$ to a polynomial of the form $p_0(h)$, where
$p_0(h) \in \bbz[\tau,\tau^{-1}][h] \subset
\bbz[\xi,\xi^{-1}][h].$ Suppose that $p(h,\scx_{2m-1}) \in \ker
(\varphi\circ \Psi_{2m-1}) $. It follows that $ \wPsi_{2m-1}(
p_0(h)) = p_0(h) \in \calc_{2m-1} $ (see~\eqref{eq:Psi-hat}) and
we can write
$$
p_0(h) = \alpha(h,\phi) (h^m - 2\phi) + \beta(h,\phi) \phi^2,
$$
where $\alpha(h,\phi), \beta(h,\phi) \in
\bbz[\xi,\xi^{-1}][h,\phi].$ Passing to
$\bbq(\xi,\xi^{-1})[h,\phi]$ and evaluating at $\phi
=\frac{1}{2}h^m$ one obtains the following equality in
$\bbq(\xi,\xi^{-1})[h]$:
$$
p_0(h) = \beta(h,\frac{1}{2}h^m)\frac{1}{4}h^{2m}.
$$
One concludes that $h^{2m} | p_0(h)$. This shows that $p_0(h) \in
\la \ve \ra + J_{2m-1}$ and hence, so does $p(h,\scx_{2m-1}).$

The case $n=2m$ is similar.
\end{proof}

In order to prove the surjectivity of $\oPsi_n$ we resort to the
spectral sequence whose $E_2$-term in described explicitly in
Remark \ref{rem:ss}.

By construction, the elements of $\cala$ are infinite cycles.
Under the description of Remark~\ref{rem:ss}, the class $h\in
CH^1(\calq_{n,\bbc})$ represents an element in $E_2^{0,2}(0)$.
Since $h$ is the image of  $c_1(\calo(1))\in
H^{2,1}(\calq_{n,\bbc};\uz)$ under the forgetful map $\varphi$, it
follows by Proposition~\ref{prop:facts}.{\bf \ref{prop:facts:iv}}
that $h$ is an infinite cycle.

The discussion above leads us to define the following subsets
 of the $E_2$-term $\cale:$
\begin{align*}
\cals_n & \equdef \{\ve^a\tau^b h^c\phi \mid a > 0, c\geq 0 \}, \quad \text{ if } n\not\equiv 0 \mod{4} \\
\cals_n & \equdef \{\ve^a\tau^b h^c\phi \mid a > 0, c>0 \}, \quad
\text{ if } n \equiv 0 \mod{4} .
\end{align*}

\begin{lemma}
\label{lem:cycles}
\begin{enumerate}[i.]
\item For each $n$, the set $\cals_n$ contains no non-zero infinite cycles.
\item If $n=2m$ , then $\chi \equdef h^{2m}-2\phi$ is an infinite cycle and
$E^{*,*}_\infty(*)$ is generated by \ $h$ and $\chi$ as an
$\cala$-algebra.
\item If $n=2m-1$, then $E^{*,*}_\infty(*)$ is generated by $h$ as an $\cala$-algebra.
\end{enumerate}
\end{lemma}
\begin{proof}
{\it i.} Suppose that $ \ve^a\tau^b h^c\phi \in \cals_n$ is an
infinite cycle. Since $h$ and $\ve$ are  infinite cycles,  this
would imply that $\ve^k\tau^b h^n \phi $ is an infinite cycle for
all $k\geq a.$ These elements live in the top row of the spectral
sequence and hence they would survive to $E_\infty,$ contradicting
the fact that $H^{r,s}_\text{Bor}(\calq_n(\bbc);\uz) = 0$ for all
$r\geq 2n;$ cf. Properties \ref{prop:borel}(ii).
\medskip

{\it ii.} Suppose now $n=2m$. Under the isomorphism
\eqref{eq:chow} the forgetful map $\varphi$ sends the element
$j_!\bone \in H^{n,-1}(\calq_{2m}(\bbc);\uz)$ to $\xi^{-1-m}(h^m
-2\phi)$. It follows by Proposition~\ref{prop:facts}(iv) that the
corresponding element $\xi^{-1-m} \chi\in \cale $ (under the
description of Remark~\ref{rem:ss}) survives to $E_\infty.$

{\it iii.} Combining the first assertion in this Lemma with Remark
\ref{rem:ss} completes the proof.
\end{proof}

It follows from the lemma that the map $\Psi_n$ is onto. Indeed,
consider first the case $n$  \emph{odd}. By the lemma,
$E^{*,*}_\infty(*)$ is generated as an $\cala$-algebra by the
image $h$ of $c_1(\calo(1))$ under the projection
$$H^{*,*}(\calq_n(\bbc);\uz)\to
Gr^0(H^{*,*}(\calq_{2m}(\bbc);\uz)).$$ Since
$c_1(\calo(1))=\Psi_n(h)$ this gives
$$
\frac{Im(\Psi_n)\cap F^j}{Im(\Psi_n)\cap F^{j-1}} =
Gr^jH^{*,*}(\calq_n(\bbc);\uz),
$$
hence $\Psi_n$ is onto. If $n=2m$ is \emph{even}, then
$E^{*,*}(*)$ is generated as an $\cala$-algebra by $h$ and the
image $x$ of $j_!\bone$ under the projection
$H^{2m,-1}(\calq_{2m}(\bbc);\uz)\to
Gr^0(H^{*,*}(\calq_{2m}(\bbc);\uz)).$ Since $j_!\bone,\,
c_1(\calo(1))\in Im(\Psi_{2m})$ it follows, as in the previous
case, that $\Psi_{2m}$ is onto.

The Theorem is now proven.

\end{proof}

\begin{corollary}
\label{cor:free}
The map $\Psi_{2m-1}$ induces an isomorphism
$$
H^{*,*}(\calq_{2m-1}(\bbc);\uz)/\text{tor}\  \cong\
\bbz[\tau,\tau^{-1},h]/\la h^{2m} \ra.
$$
Similarly, $\Psi_{2m}$ induces an isomorphism
$$
H^{*,*}(\calq_{2m }(\bbc);\uz)/\text{tor} \cong
\bbz[\tau,\tau^{-1},h,\chi ]/\la h^{2m+1}, h\chi ,
\tau^{m+1}\chi^2-(-1)^m h^{2m} \ra
$$
\end{corollary}
\begin{proof}
Left to the reader.
\end{proof}

%%%%%%%%%%%%%%%%%%%%%%%%%%%%%%%%%%%%%%%%%%%%%%%%%%%%%%%%%%%%%%%%%%%%%%%%%%%%%%%%%%

%\end{document}

%\input{newisotropic.tex}

\section{Isotropic quadrics}
\label{sec:isotropic}

Consider a real quadratic form $\mbq$ of rank $n+2$ that can be
written as $\mbq=\mbq'+\mbh,$ where $\mbh$ is a hyperbolic factor.
Let $P \in X_\mbq(\bbr)$ be a real point in $X_\mbq$ and denote by
$T_PX_\mbq $ the projective tangent space to $X_\mbq$ at $P.$ The
intersection $X_\mbq\cap T_PX_\mbq = \cxs X_{\mbq'}$ is the ruled
join of the quadric $X_{\mbq'}$ and the point $P.$ In particular,
the set of complex points $\cxs X_{\mbq'}(\bbc)$ in the analytic
topology is the Thom space of the Real bundle
$\calo_{X_{\mbq'}}(1)(\bbc) $ over $X_{\mbq'}(\bbc),$ and $X_\mbq
- \cxs X_{\mbq'}$ is isomorphic to affine space $\bba^n.$ Let $P_0
\in X_\mbq-\cxs X_{\mbq'}$ correspond to $0 \in \bba^n$ under this
isomorphism.

\begin{definition}
\label{def:filtrations}
Denote
\begin{equation*}
\begin{array}{lll}
Y_2   = \{ P_0 \}, & Y_1 = X_{\mbq'}, & Y_0 = \{ P \},
\quad \text{ and } \\
X_2  = X_\mbq, & X_1  = \cxs X_{\mbq'}, &  X_0 = \{ P \} .
\end{array}
\end{equation*}
Note that the filtration $X_\mbq=X_2 \supset X_1  \supset X_0
\supset X_{-1} =\emptyset$ comes with projections $\pi_i \colon
X_i - X_{i-1} \to Y_i$ that are real algebraic vector bundles, and
satisfies  $a_2 :=\operatorname{codim}{X_2} = 0, \quad a_1
:=\operatorname{codim}{X_1} =  1$ and $a_0
:=\operatorname{codim}{X_0} = n,$ respectively.
\end{definition}

\begin{notation}
Given a bigraded $\calb$-module $N = \oplus_{r,s} N^{r,s}$ and
integer $a$, denote by $N(a)$ the bigraded $\calb$-module whose
summand of degree $(r,s)$ is $N^{r+2a,s+a}.$
\end{notation}

\begin{proposition}
\label{prop:isotropic}
Given a quadratic form $\mbq = \mbq'+\mbh$ as above,   denote
$\bbh \equdef H^{*,*}(X_\mbq(\bbc);\uz)$ and $\bbh' \equdef
H^{*,*}(X_{\mbq'}(\bbc);\uz)$. Then the Bredon cohomology of
$X_\mbq(\bbc)$ has a decomposition
$$
\bbh \cong \calb \oplus \bbh'(-1) \oplus \calb(-n),
$$
as a sum of bigraded $\calb$-modules.
\end{proposition}
\begin{proof}
Denote $H^{(p)}(X) \equdef \oplus_{r\in \bbz}
H^{r,p}(X(\bbc);\uz)$, whenever $X$ is a real variety and $p\in
\bbz$. The cycle map \eqref{eq:cycle_map} $CH^p(X) \to H^{(p)}(X)$
gives $H^{(*)}(X) = \oplus_p H^{(p)}(X)$ the structure of a
$CH^*(X)$-algebra, whenever $X$ is smooth. It is easy to see that,
in the category $\calv_\bbr$ of smooth proper real varieties, the
functor $H^{(*)}$ becomes a graded geometric cohomology theory in
the sense of \cite{Kar-flag}.  The proposition now follows from
Definition \ref{def:filtrations}, Theorem 6.5 and Corollary 6.11
of \cite{Kar-flag}.
\end{proof}

In order to understand the product structure in the Bredon
cohomology, we provide a brief description of the isomorphism in
the proposition above.

Define $U_i \subset Y_i \times X_\mbq$  as $U_i = \{ (y,x) \mid
x\in X_i - X_{i-1} \text{ and } \pi_i(x) = y \}, $ cf. Definition
\ref{def:filtrations}, and let $\Gamma_i = \overline{U_i}$ denote
its closure. Then, $\Gamma_i$ defines a correspondence $ \Gamma_i
\in CH^{\dim{Y_i}+a_i}(Y_i \times X_\mbq).$ Each such
correspondence defines a map of (bigraded) $\calb$-modules
$$
\gamma_i \colon H^{*,*}(Y_i(\bbc);\uz)(-a_i) \to
H^{*,*}(X_\mbq(\bbc);\uz)
$$
and the isomorphism above is simply given by the  sum $\gamma_2 +
\gamma_1 + \gamma_0.$ Recall that $\gamma_i$ is explicitly defined
on a cohomology class $\alpha \in H^{*,*}(Y_i(\bbc);\uz)$ as
follows. Let $Y_i \xleftarrow{pr_1} Y_i\times X_\mbq
\xrightarrow{pr_2} X_\mbq$ be the projections. Then, invoking
Poincar\'e duality, $\gamma_i $ is uniquely determined by
$$
\gamma_i(\alpha) \cap [X_\mbq(\bbc)] = pr_{2 *} \left(
pr_1^*\alpha \cap [\Gamma_i(\bbc)] \right),
$$
where $[ - ] $ denotes the fundamental class of a subvariety in
Bredon homology.

These constructions hold in great generality, whenever the variety
$X_\mbq$ has a filtration satisfying the conditions in Definition
\ref{def:filtrations}. In our case, $X_2 - X_1 \cong \bba^n $ and
the projection $\pi_2 \colon X_2-X_1 \to \{ P_0 \}$ is identified
with the ``retraction'' from $\bba^n$ to $0 \in \bba^n.$ Hence,
$\Gamma_2 = \{ pt \} \times X_\mbq$ and the map $\gamma_2 \colon
\calb \to H^{*,*}(X_\mbq(\bbc);\uz)$ splits off the cohomology of
a point via the real point $P_0 \in X_\mbq(\bbr).$ Similarly, $
X_0 - X_{-1} = \{ P \}$ and $\Gamma_0 = \{ (P ,P ) \} \subset \{ P
\} \times X_\mbq$, and the map $\gamma_0 \colon \calb(-n) \to
H^{*,*}(X_\mbq(\bbc);\uz)$ is the Gysin map associated to the
inclusion $\{ P \} \hookrightarrow X_\mbq.$

Finally, $X_1- X_0$ is the total space of the hyperplane bundle
$\calo_{X_{\mbq'}}(1)$ over $X_{\mbq'}.$ It is easy to see that
$\Gamma_1 \subset X_{\mbq'}\times X_\mbq$ is  the projective
closure $\bbp(\calo_{X_{\mbq'}}(1) \oplus \bone )$ of
$\calo_{X_{\mbq'}}(1).$ The restriction of the projection $pr_2$
to $\Gamma_1(\bbc) $ is the composition $\Gamma_1(\bbc) \to
Thom(\calo(1)(\bbc)) = X_1(\bbc) \hookrightarrow X_\mbq(\bbc)$,
where the first map corresponds to collapsing the section at
infinity and the last one is just the inclusion $i_1 \colon
X_1(\bbc) \hookrightarrow X_\mbq(\bbc).$

\begin{remark}
It is easy to see that $\gamma_1 (\alpha ) \cap [X_\mbq] = i_{1*}
( \calt( \alpha \cap [X_{\mbq'}] )),$ where $\calt \colon
H_{*,*}(X_{\mbq'}(\bbc);\uz) \to
H_{*+2,*+1}(Thom(\calo(1)(\bbc));\uz)$ is the Thom isomorphism in
Bredon homology.
\end{remark}
\begin{proposition}
\label{prop:facts-chern}
Using the notation above, one has $$\gamma_1 (\bone) = \uc_1\left(
\calo_{X_\mbq}(1)(\bbc)\right) \in H^{2,1}(X_\mbq(\bbc);\uz).$$
\end{proposition}
\begin{proof}
By definition $\gamma_1(\bone) \cap [X_\mbq] = pr_{2*}\left(
pr_1^* \bone \cap [\Gamma_1] \right) =  pr_{2*} [\Gamma_1] . $ The
latter element is precisely $[pr_2(\Gamma_1)] = [X_1]$, a
hyperplane section of $X_\mbq$.
\end{proof}

Consider the real quadric $\calq_{n,s},$ and assume $2s \leq n .$
In order to simplify notation, denote the Bredon cohomology ring
of $\calq_{n,s}(\bbc)$ by
\begin{equation}
\label{eq:cohomology}
\bbh_{n,s} \equdef H^{*,*}(\calq_{n,s}(\bbc);\uz).
\end{equation}
Using Propositions \ref{prop:isotropic} and
\ref{prop:facts-chern}, together with an induction argument one
easily obtains the additive structure of $\bbh_{n,s}$.

\begin{proposition}
\label{prop:explicit}
Let $\mbh \in H^{2,1}(\calq_{r,s}(\bbc);\uz)$ denote the first
Chern class $\uc_1(\calo(1))$ of $\calo_{\calq_{n,s}}(1)$ and let
${\eta} \in H^{2(n+1-s),n+1-s}(\bbq_{n,s}(\bbc);\uz)$ be the
Poincar\'e dual to a $\bbp_\bbr^{s-1} \subset \calq_{n,s} $. Then
one has an inclusion of $\calb$-modules
\begin{equation}
\label{eq:dagger}
\jd \colon \bbh_{n-2s,0} \to \bbh_{n,s}(s)
\end{equation}
which induces an identification of bigraded $\calb$-modules:
\begin{align*}
\bbh_{n,s}\ = \ & \  \calb \cdot \bone \oplus \calb \cdot \mbh
\oplus \cdots \oplus \calb \cdot \mbh^{s-1}\ \oplus \
\bbh_{n-2s,0}(-s) \ \oplus \\ & \ \calb \cdot \eta \oplus \calb
\cdot\mbh \eta \oplus \cdots \oplus \calb\cdot \mbh^{s-1} \eta.
\end{align*}
\end{proposition}

\subsection{The multiplicative structure}

We proceed to determine the multiplicative structure of the
cohomology ring of $\calq_{n,s}$, with $n\geq 2s.$ The terms
\emph{vector bundles} and \emph{spaces} henceforth mean \emph{Real
vector bundles} and \emph{Real spaces}; cf. \cite{Ati66}.

Given $0\leq s \in \bbz$, define a subring of the polynomial ring
$\calb[h,\scx,y]$ by
\begin{equation}
\label{eq:def_Bs-1}
\calb_s[h,\scx,y] \equdef \calb[h]+ \la h^s \ra,
\end{equation}
where $\la h^s \ra$ is the ideal generated by $h^s$. Every $P \in
\calb_s[h,\scx,y]$ can be written uniquely as
\begin{equation}
\label{eq:poly}
P = P_0 + h^s \cdot P_1,
\end{equation}
with $ P_0 =a_0\cdot \bone + a_1\cdot h + \cdots +
a_{s-1}h^{s-1},\ a_i \in \calb$, $i=1,\ldots,s-1$ and $P_1 \in
\calb[h,\scx,y].$

In what follows, we identify  $\bbh_{n-2s,0}$ with
$\cala[h,\scx]/J_{n-2s}$ and use the surjection of
$\calb$-algebras $\calb[y]\to \cala$ of $\calb$-algebras; cf.
Theorem \ref{thm:MAIN} and Remark \ref{rem:fg}(ii), respectively.
Let $\Lambda(\eta)$ be the exterior algebra over $\bbz$ on one
generator $\eta$. It follows directly from Proposition
\ref{prop:explicit} and \eqref{eq:poly}that the map of
$\calb$-\emph{modules} defined by
\begin{align}
\label{eq:themap}
\Psi   \colon \calb_s[h,\scx,y]\otimes \Lambda(\eta) & \longrightarrow\ \bbh_{n,s} \notag \\
h^j\otimes \eta^\epsilon &\longmapsto\ \mbh^j\eta^\epsilon, \quad j=0,\ldots, s-1, \quad \epsilon= 0,1 \notag \\
h^s P \otimes \eta & \longmapsto\ 0, \quad \text{ for all } P \in \calb[h,\scx, y] \notag \\
h^{s+r}\scx^k y ^a \otimes 1 & \longmapsto\ 0 , \quad \text{ if
either } kr\neq 0, \text{ or } k>1 \notag \\
h^{s+r}y^a \otimes 1 & \longmapsto \mbh^r \jd(\tau^{-a}), \quad \text{ for all } a \geq 0 \notag \\
h^{s}\scx  y^b \otimes 1 & \longmapsto \jd(\tau^{-b}\scx ), \quad \text{ for all } b \geq 0\notag \\
\end{align}
is surjective. In what follows, we show that $\Psi$ is indeed a
$\calb$-\emph{algebra} homomorphism and identify its kernel.
\bigskip

Let $L_1$ and $L_2$ be line bundles over an $H\uz$-oriented
compact manifold $X$ with first Chern classes $t_1, t_2 \in
H^{2,1}(X;\uz)$, respectively. Let $\pi \colon \bbp \to X$ denote
the projection from $\bbp \equdef \bbp(L_1\oplus L_2)$ onto $X$
and let $\si \colon X \to \bbp$ be the ``section at infinity'',
where $X$ is identified with $\bbp(0\oplus L_2)$. The following
identities on characteristic classes are well-known.

\begin{lemma}
\label{lem:0}
Using the notation above, let $\zeta \equdef c_1(\calo_{L_1\oplus
L_2}(1)) \in H^{2,1}(\bbp,\uz).$ Then
\begin{equation}
\label{eq:section1}
\sic\calo_{L_1\oplus L_2}(1) = \dual{L}_2,
\end{equation}
where $\dual{L}_2$ is the dual of $L_2$, and
\begin{equation}
\label{eq:section2}
\zeta \cap [\bbp] = - \pi^* t_1 \cap [\bbp] + \sih[X].
\end{equation}
Furthermore, for all $r\geq 0$ one has:
\begin{equation}
\label{eq:zeta}
\zeta^{r+1} \cap [\bbp] = (-1)^{r+1} \pi^* t_1^{r+1} \cap [\bbp] \
+\ (-1)^r \sum_{i+j=r}\ \sih (t_1^it_2^j \cap [X]).
\end{equation}
\end{lemma}
\begin{proof}
The first two identities are standard facts about characteristic
classes; cf. \cite{fulton-IT}.

The last assertion is true for $r=0;$ cf. \eqref{eq:section2}.
Assume true for $r_0-1$, then
\begin{align*}
\zeta^{r_0+1} \cap [\bbp]  & = (-1)^{r_0} \pi^*t_1^{r_0} \cap
\zeta \cap [\bbp] \\ & \quad + (-1)^{r_0-1} \sum_{i+j=r_0-1}\
\zeta\cap \sih
(t_1^it_2^j\cap [X] ) \\
& = (-1)^{r_0} \pi^*t_1^{r_0} \cap \left\{ - \pi^*(t_1) \cap
[\bbp] + \sih [X] \right\} \\
&\quad  + (-1)^{r_0-1} \sum_{i+j=r_0-1}\
\sih \left( \sic\zeta \cap (t_1^it_2^j\cap [X] ) \right) \\
& \overset{\eqref{eq:section1}}{=} (-1)^{r_0+1} \pi^*(t_1^{r_0+1})
\cap [\bbp]  \\ & \quad + (-1)^{r_0} \sih(t_1^{r_0} \cap [X]) +
(-1)^{r_0}
\sum_{i+j = r_0-1} \sih t_1^it_2^{j+1} \cap [X] \\
& =  (-1)^{r_0+1} \pi^* t_1^{r_0+1} \cap [\bbp] \ +\ (-1)^{r_0}
\sum_{i+j=r_0}\ \sih (t_1^it_2^j \cap [X]).
\end{align*}
\end{proof}
\bigskip

Denote $\calq = \calq_{n,s}$ and $\calq'=\calq_{n-2s,0}$, and
consider the following diagram
$$
 \xymatrix{ & \bbp(\pi_1^* \calo(1) \oplus \pi_2^*\calo(1) )
\ar[d]_{\pi} \ar[ddrr]_{b'} \ar[ddrrrr]^{b} & & & & \\
& \calq' \times \bbp^{s-1}_\bbr \ar[dl]^{\pi_1} \ar[dr]_{\pi_2} & & & & \\
\calq' & & \bbp^{s-1}_\bbr & \bbp^{s-1}_\bbr \# \calq' \ar[rr]_j &
& \calq , }
$$
where $\bbp^{s-1}_\bbr \# \calq'$ is the ruled join of
$\bbp^{s-1}_\bbr$ and $\calq'$, $j$ is the natural inclusion, $b'$
is the blow-up map of $\bbp^{s-1}_\bbr \# \calq'$ along
$\bbp^{s-1}_\bbr \amalg \calq'$ and $b = j\circ b'.$

\begin{proposition}
\label{def:jd}
The map
\begin{equation*}
\jd \colon \bbh_{n-2s,0} \longrightarrow \bbh_{n,s}(s)
\end{equation*}
introduced in \eqref{eq:dagger}  sends $\alpha \in
H^{p,q}(\calq')$ to $\jd(\alpha)$ satisfying
$$
\jd(\alpha)\cap [\calq] \equdef (-1)^q b_*(\pi^*(\alpha \times
\bone) \cap [\bbp] ) \ \in H^{p+2s,q+s}(\calq),
$$
where $\bbp \equdef \bbp(\pi_1^* \calo(1) \oplus \pi_2^*\calo(1))
.$ In particular, $\jd(1) = \mbh^s.$
\end{proposition}
\begin{proof}
Left to the reader.
\end{proof}

\begin{notation}
\label{not:proj}
For $X = \calq'\times \bbp^{s-1}_\bbr,$  let $L_1=\pi_1^*\calo(1)$
and $L_2= \pi_2^*\calo(1)$, and denote $c_1L_1 = t_1 = h \times
\bone$ and $c_1L_2 = t_2 = \bone \times t.$ Let  $i \colon
\bbp^{s-1}_\bbr \hookrightarrow \calq$ be the inclusion.
\end{notation}

\begin{lemma}
\label{lem:1}
Given $\alpha \in H^{p,q}(\calq'(\bbc);\uz)$ and $r\geq 0$ one
has:
\begin{multline*}
\mbh^{r+1}\cap \jd(\alpha)\cap [\calq(\bbc)] = \\
\jd(h^{r+1}\alpha)\cap [\calq(\bbc)] + (-1)^{r+q} i_* \pi_{2*}
\sum_{i+j = r}\ (\alpha h^i \cap [\calq'(\bbc)])\times (t^j \cap
[\bbp^{s-1}_\bbr(\bbc)]).
\end{multline*}
\end{lemma}

\begin{proof}
The fact that $b^*\mbh = \zeta$ and Lemma \ref{lem:0} give:
\begin{align*}
\mbh^{r+1} & \cap \jd(\alpha)\cap [\calq(\bbc)] = (-1)^q
\mbh^{r+1}\cap
b_*( \pi^*(\alpha\times \bone )\cap [\bbp(\bbc)] )  \\
& = (-1)^q b_*(
\pi^*(\alpha\times \bone )\cap \{ b^*\mbh^{r+1} \cap [\bbp(\bbc)]\}) \\
& =(-1)^q b_* ( \pi^*(\alpha\times \bone )\cap  \{ (-1)^{r+1}
\pi^*(h^{r+1} \times \bone) \cap [\bbp(\bbc)] \\
& + (-1)^r \sih \sum_{i+j=r} (h^i\times t^j)\cap [(\calq'\times
\bbp^{s-1}_\bbr)(\bbc)] \ \} \ )  \\
& =(-1)^{q+r+1} b_*\left( \pi^*(h^{r+1}\alpha \times \bone ) \cap
[\bbp(\bbc)] \right)  \\
& + (-1)^{r+q} b_*\left( \pi^*(\alpha \times \bone ) \cap \sih
\sum_{i+j=r} (h^i\times t^j)\cap [ (\calq' \times
\bbp^{s-1}_\bbr)(\bbc) ] \right)
\\
 & = \jd( h^{r+1}\alpha ) \cap [\calq]  \\
 & \ + (-1)^{r+q} b_*\sih \left(
\sum_{i+j=r} (\alpha h^i\times t^j)\cap [(\calq'\times
\bbp^{s-1}_\bbr)(\bbc)] \right)
\\
& = \jd( h^{r+1}\alpha ) \cap [\calq(\bbc)]  \\
&\ + (-1)^{r+q} i_* \pi_{2*} \left( \sum_{i+j=r} (\alpha h^i\times
t^j)\cap [(\calq'\times \bbp^{s-1}_\bbr)(\bbc)] \right) .
\end{align*}
\end{proof}

\begin{corollary}
\label{cor:hTa}
For all $a\geq 0$ one has:
\begin{equation*}
\mbh^r \jd(\tau^{-a} h^{r'}) =
\begin{cases}
\jd(\tau^{-a}h^{r+r'}) &, \text{ if } 0 \leq r+r' \leq n-2s \\
2\tau^{-a} \eta \mbh^{r+r'-1 - n + 2s} &, \text{ if } n-2s +1 \leq
r + r' \leq n-s \\
0 &, \text{ if } n-s < r+r'.
\end{cases}
\end{equation*}
\end{corollary}
\begin{proof}
One has
\begin{align*}
\mbh^r\jd( & \tau^{-a}h^{r'}) \cap [\calq(\bbc)]  \\ & =
\jd(\tau^{-a} h^{r+r'} )\cap [\calq(\bbc)]  \\
& + (-1)^{r+r'-1}
\sum_{i+j=r-1} i_*\pi_{2*}\left(
\tau^{-a}h^{i+r'}\cap[\calq'(\bbc)]\ \times \ t^j \cap
[\bbp^{s-1}_\bbr](\bbc) \right) .
\end{align*}
For $r+r'-1< n-2s$ the terms in the summation above vanish for
dimensional reasons, hence $ \mbh^r\jd( \tau^{-a} h^{r'})  =
\jd(\tau^{-a}h^{r+r'}) $ in this case.
For $r+r'-1\geq n-2s$ one has $h^{r+r'}= 0$ in $\bbh_{n-2s,0}$,
thus giving
\begin{multline}
\mbh^r\jd(\tau^{-a}h^{r'}) \\  = (-1)^{r-1} i_*\pi_{2*} \left(
\tau^{-a}h^{n-2s} \cap [\calq'(\bbc)] \times t^{r+r'-1-n+2s} \cap
[\bbp^{s-1}_\bbr(\bbc)] \right).
\end{multline}
Since $\dim \bbp^{s-1}_\bbr =
s-1$, this immediately shows that $\mbh^r\jd(\tau^{-a}h^{r'}) =
0,$ if $r+r'>n-s.$

Now, if $n-2s+1\leq r+r' \leq n-s$, i.e. $0\leq r+r'-1-n-2s \leq
s-1,$ one has
$$
i_* (t^{r+r'-1-n+2s} \cap \bbp^{s-1}_\bbr(\bbc)) = \eta
\mbh^{r+r'-1-n+2s} \cap [\calq(\bbc)].
$$

On the other hand, for all $m\geq 0$ one has
$$\pi_{2*}(\tau^{-a}h^{n-2s}\cap [\calq'(\bbc)] \times t^m\cap
[\bbp^{s-1}_\bbr(\bbc) ] ) = (2\tau^{-a}) t^m\cap
[\bbp^{s-1}_\bbr(\bbc)]. $$
Therefore,
$$
\mbh^r\jd(\tau^{-a}h^{r'}) = (2\tau^{-a})\eta \mbh^{r+r'-1-n+2s},
$$
for $n-2s+1\leq r+r' \leq n-s.$
\end{proof}
\begin{remark}
Note that $2\tau^{-a}$ is an element in $\calb$ for all $a\in
\bbz$ while $\calb$ does not have an inverse $\tau^{-a}$ for
$\tau^a$, for $a> 0$.
\end{remark}

\begin{lemma} For all $\alpha \in \bbh_{n-2s,0}$ and $k \in \bbz$ one has
\label{lem:2}
\begin{enumerate}[a.]
\item $\jd(\alpha) \cdot \eta = 0$;
\item $b^*(\jd (\tau^k \alpha) ) = \tau^k b^*(\jd(\alpha)) \in H^{*,*}(\bbp(\bbc),\uz).$
\end{enumerate}
\end{lemma}
\begin{proof}
One has
\begin{multline}
\eta \jd(\alpha) \cap [\calq(\bbc)]\\ = \pm \ \eta \cap b_*(
\pi^*(\alpha \times \bone) \cap [\bbp(\bbc)])  = \pm\
b_*(\pi^*(a\times \bone) \cap b^*\eta \cap [\bbp(\bbc)]).
\end{multline}
However, one has $\deg{\eta} = ( 2(n-2+1),n-2+1)$, while
$\dim{\bbp} = n-s,$ and hence $b^*\eta = 0,$ thus proving the
first assertion of the Lemma.

Now, observe that the second assertion holds for $k\geq 0,$ since
$\jd$ is a homomorphism of $\calb$-modules. If $k = -a,$ with
$a>0$, one has $ \tau^a \jd(\tau^{-a} \alpha) = \jd(\alpha), $ and
hence $b^*(\tau^a \jd(\tau^{-a} \alpha)) = b^*\jd(\alpha).$ This
is equivalent to saying that $\tau^a b^*(\jd(\tau^{-a}\alpha)) =
b^* \jd(\alpha).$ Since $\symg{2}$ acts freely on $\bbp(\bbc),$
$\tau $ is invertible in the cohomology of $\bbp(\bbc),$ and one
has $b^*(\jd(\tau^{-a}\alpha)) = \tau^{-a}b^*(\jd\alpha).$
\end{proof}

\begin{corollary}
\label{lem:4}
Denote $\scx_a \equdef \jd(\tau^{-a}\scx)$ and $\sct_a
=\jd(\tau^{-a})$. Then:
\begin{enumerate}[a.]
\item $\sct_a\sct_b = \sct_{a+b}\mbh^s$,
\item $\mbh^{r+1}\scx_a = 0,$ for all $r\geq 0,$
\item $\sct_a \scx_b = 0$ for all $a,b \geq 0,$
\item $\scx_a\scx_b = 0$ for all $a, b \geq 0$.
\end{enumerate}
\end{corollary}
\begin{proof}
To prove {\it \textbf{a}} one observes that
\begin{align*}
\sct_a\sct_b \cap [ \calq(\bbc)] & = \sct_a \cap
b_*(\pi^*(\tau^{-b}\times \bone)\cap [\bbp(\bbc)]) \\
& =  b_*(\tau^{-a}b^*(\jd(\bone) ) \pi^*(\tau^{-b} \times
\bone)\cap
[\bbp(\bbc)]) \\
& = b_*(b^*\mbh^s \cdot  \pi^*(\tau^{-a-b} \times \bone)\cap [\bbp(\bbc)]) \\
& = \mbh^s \sct_{a+b} \cap [\calq(\bbc)].
\end{align*}

To prove assertion {\it \textbf{b}} one has:
\begin{align*}
\mbh^{r+1} \scx_a & \cap [ \calq(\bbc)]   =
 \pm b_*(\ \pi^*(\tau^{-a}\scx \times \bone) \cap\{ b^*\mbh^{r+1}
\cap [\bbp(\bbc)] \} \ ) \\
 & = \pm b_*  ( \pi^* ( \tau^{-a}\scx \times \bone) \cap \{
(-1)^{r+1} \pi^*(h^{r+1} \times \bone ) \cap [\bbp(\bbc)]  \\
&   \quad + (-1)^{r} \sih
\sum_{i+j=r} (h^i\times t^j)\cap [(\calq'\times \bbp^{s-1}_\bbr)(\bbc)] \ \} \ ) \\
 & = \pm b_* ( \pi^* ( \tau^{-a}h^{r+1} \scx \times \bone) \cap
[\bbp(\bbc)]  \\
& \quad \pm  \sih \sum_{i+j=r} (\tau^{-a}\scx h^i \times t^j)\cap
[(\calq'\times \bbp^{s-1}_\bbr)(\bbc)]  )  \\
  & =  \pm i_* \pi_{2*} \left\{ ( \tau^{-a} \scx \times t^r ) \cap
[(\calq' \times \bbp^{s-1}_\bbr)(\bbc)] \right\} \\
  & = \pm i_* \pi_{2*} \left\{  \tau^{-a}\cap [S^{n-2s}_a] \times t^r
\cap [\bbp^{s-1}_\bbr(\bbc)]  \right\}.
\end{align*}
The latter expression is clearly zero if $n>2s$ (K\"unneth formula
for  $\calq' \times \bbp^{s-1}_\bbr$), and when $n=2s$ one has
$S^0_a \cong \symg{2} = \{ 0, 1 \}$ and $[S^0_a] = [1]-[0] $. It
follows that $\pi_{2*}(\tau^{-a}\cap [S^0_a] \times t^r \cap
[\bbp^{s-1}_\bbr] ) = 0,$ as well.

To prove assertion {\it \textbf{c}} one has
\begin{align*}
\sct_a\scx_b \cap [\calq(\bbc)]   & = \pm \sct_a \cap
b_*(\pi^* ( \tau^{-b} \scx \times \bone ) \cap [\bbp(\bbc)] ) \\
&  = \pm b_*( b^*\jd(\tau^{-a}) \cap \pi^* ( \tau^{-b} \scx
\times \bone ) \cap [\bbp(\bbc)] ) \\
 & = \pm  b_*(\tau^{-a} b^*\mbh^s \cap \pi^* ( \tau^{-b} \scx
\times \bone ) \cap [\bbp(\bbc)] ) \\
 & = \pm  b_*(\pi^*(\tau^{-a-b}h^s \scx \times \bone ) \cap
[\bbp(\bbc)] ) \\
\ & \quad \pm i_* \pi_{2*} \left( \sum_{i+j=s-1} \tau^{-a-b} \scx
h^i \cap
[\calq'(\bbc)] \times t^j \cap [\bbp^{s-1}_\bbr(\bbc)] \right) \\
&  = \pm i_* \pi_{2*} \left( \tau^{-a-b} \scx \cap [\calq'(\bbc)]
\times t^{s-1} \cap [\bbp^{s-1}_\bbr(\bbc)] \right) = 0,
\end{align*}
using the arguments of the previous lemma.

Observe that $\symg{2}$ acts freely on $\bbp(\bbc)$ and that
$\dim{\bbp} = n-s$. Hence, whenever $\alpha \in
H^{i,j}(\bbp(\bbc);\uz)$ and $i> 2(n-s)$ one has $\alpha =0$; cf.
Property \ref{prop:borel}(ii). Since
\begin{align*}
\scx_a\cap ( \scx_b \cap [\calq(\bbc)]) & = \pm \scx_a \cap
b_*(\pi^* (\tau^{-b} \scx \times 1) \cap [\bbp(\bbc)]) \\
& = \pm b_*( \{ b^*\scx_a \cdot
 \pi^* (\tau^{-b} \scx \times 1 ) \} \cap [\bbp(\bbc)] )
\end{align*}
and $\deg{b^*\scx_a \cdot
 \pi^* (\tau^{-b} \scx \times 1 )} = (2n, 2(s-1)-a-b) $,
 it follows that $\scx_a\cdot \scx_b = 0,$ thus showing the last
 assertion of the lemma.
\end{proof}

\begin{proposition}
\label{prop:final}
Write $n=2m-\delta.$ Then
\begin{enumerate}[i.]
\item The map $\Psi \colon \calb_s[h,x,y]\otimes \Lambda(\eta) \to
H^{*,*}(\calq_{n,s}(\bbc);\uz)$ \eqref{eq:themap} is a surjective
homomorphism of $\calb$-algebras.
\item The kernel of $\Psi$ is the ideal
$$\ \cali_{n,s} = [h^s]\cdot \tilde{J}_{n-2s} \otimes \Lambda(\eta)\ +\ [h^s]\otimes
\la \eta \ra  \ + \ \la h^{n-s+1}\otimes 1 - 2(1\otimes \eta) \ra
,\
$$
where $\tilde{J}_{n-2s} = [g_1,g_2, g_3,g_4,g_5] \subset
\calb[h,\scx,y]$ is the ideal generated by the elements
\begin{align*}
g_1 & = \mbf_{m-s}, \ \ g_2 = \ve^{1-\delta} \tau^{m-s} \scx -
h^{1-\delta} \mbf_{m-s-1}, \ \ g_3 = h \scx, \\
g_4& = h^{2(m-s)} - \delta \{ (-1)^{m-s} \tau^{m-s+1} \scx^2 \} ,
\ \ \text{ and } \ \ g_5 = \tau y - 1 .
\end{align*}
\end{enumerate}
\end{proposition}
\begin{proof}
The fact that $\Psi$ is a ring homomorphism follows directly from
Corollaries \ref{cor:hTa}, \ref{lem:4} and equation
\eqref{eq:poly}, together with the obvious fact that $\eta^2 = 0$
in $\bbh_{n,s}$. This concludes the proof of {\it \textbf{i}},
since $\Psi$ was shown to be surjective in \eqref{eq:themap}.

Let $\varrho \colon \calb[h,\scx,y] \to \bbh_{n-2s,0}$ denote the
surjection induced by  $\calb[y] \to \cala$ and the presentation
$\bbh_{n-2s,0}\cong \cala[h,\scx]/J_{n-2s}.$ It follows from
Theorem \ref{thm:MAIN} and Remark \ref{rem:fg}(ii) that the ideal
$\tilde{J}_{n-2s}$ in the statement of the proposition is
precisely the kernel of $\varrho$.

It follows from the definition of $\Psi$ and Corollary
\ref{cor:hTa} that whenever the highest power of $h$ in $P\in
\calb[h,x,y]$ (denoted $\deg_h{P}$) is less or equal than $n-2s$
then
\begin{equation}
\label{eq:final-aux}
\Psi( h^s P \otimes \ 1 ) = \jd(\varrho(P)).
\end{equation}

Now, the fact that $\deg_h{\mbf_r}\leq r$, for all $r$, together
with \eqref{eq:final-aux} and the definition of $\Psi$ shows that
$[h^s]\cdot \tilde{J}_{n-2s}\otimes \Lambda(\eta) \subset
\ker{\Psi}.$ Also, Corollary \ref{cor:hTa} gives $\mbh^{n-s+1} =
\mbh^{n-2s+1}\jd(1) = 2\eta$ and hence, the element
$h^{n-s+1}\otimes 1 - 2(1\otimes \eta)$ belongs to the kernel of
$\Psi.$  Since $[h^s]\otimes \la \eta \ra \subset \ker{\Psi}$, by
definition, one concludes that $\cali_{n,s} \subset \ker{\Psi}.$

Consider $\mbu \in \ker{\Psi}.$ Since  $[ \mbh^s ] \otimes \la
\eta \ra \subset \cali_{n,s}$ then
$$\mbu \equiv A_0\otimes 1 + B_0\otimes \eta + h^s P \otimes 1  \ \mod{\cali_{n,s}}$$
where $A_0=a_0 + a_1 h + \cdots a_{s-1}h^{s-1},\ B_0=b_0 + b_1 h +
\cdots b_{s-1}h^{s-1},$ and $P$ is an arbitrary element of
$\calb[h,x,y].$

Let us now write $P = P_0 + h^{n-2s+1}P_1 + h^{n-s+1}P_2,$ where
$\deg_h{P_0}\leq n-2s$ and $\deg_hP_1\leq s-1$. A simple
inspection shows that that $h^{n-s+1} \in \tilde{J}_{n-2s}$.
Therefore, $h^sP\otimes 1 \equiv h^sP_0\otimes 1 +
h^{n-s+1}P_1\otimes 1 \ \mod{\cali_{n,s}}.$ On the other hand,
$h^{n-s+1}P_1\otimes 1 = (h^{n-s+1}\otimes 1)(P_1\otimes 1) \equiv
2(1\otimes \eta)(P_1\otimes 1) \mod{\cali_{n,s}}.$ It follows that
one can finally write
\begin{equation}
\label{eq:finalform}
\mbu \cong A_0\otimes 1 + B'_0\otimes \eta + h^sP_0\otimes 1 \
\mod{\cali_{n,s}},
\end{equation}
where $B'_0 = B_0 + 2P_1$ and $\deg_h B'_0\leq s-1.$

Finally, it follows from \eqref{eq:final-aux} that
$$
0=\Psi(\mbu) = \Psi(A_0\otimes 1) +\Psi(B'_0\otimes \eta) +
\jd(\varrho(P_0)).
$$
Since the decomposition in Proposition \ref{prop:explicit} is a
direct sum, one concludes that $A_0 = B'_0=0$ and
$\jd(\varrho(P_0)) = 0.$ Since $\jd$ is injective, one concludes
that $P_0 \in \ker{\varrho}$, in other words $P_0 \in
\tilde{J}_{n-2s}.$ This shows that $\mbu \in \cali_{n,s}$ and
hence $\ker{\Psi} = \cali_{n,s}$
\end{proof}

With this proposition, we conclude the proof of Theorem A, stated
in the Introduction.

%%%%%%%%%%%%%%%%%%%%%%%%%%%%%%%%%%%%%%%%%%%%%%%%%%%%%%%%%%%%%%%%%%%%%%%%%%%%%%%%%%

%\input{relations.tex}

%\input{appendix.tex}

%\input{spheres.tex}

%\bibliography{quadrics-final}
\bibliographystyle{amsalpha}

\providecommand{\bysame}{\leavevmode\hbox
to3em{\hrulefill}\thinspace}
\providecommand{\MR}{\relax\ifhmode\unskip\space\fi MR }
\providecommand{\MRhref}[2]{%
  \href{http://www.ams.org/mathscinet-getitem?mr=#1}{#2}
} \providecommand{\href}[2]{#2}

\end{document}